\documentclass[10pt,notitlepage,twoside,a4paper]{amsart}
 \usepackage{amsfonts}

\usepackage{amsmath,amssymb,enumerate}

\usepackage{epsfig,fancyhdr,color}

\usepackage{amssymb}
\usepackage{amsmath,amsthm}
\usepackage{latexsym}
\usepackage{amscd}
\usepackage{psfrag}
\usepackage{graphicx}
\usepackage[latin1]{inputenc}
\usepackage[all]{xy}
\usepackage[mathcal]{eucal}

\definecolor{NoteColor}{rgb}{1,0,0}

\renewcommand{\textsc}{\textcolor{red}}

\newtheorem*{theorem 1}{\rm\bf Proposition 1}
\newtheorem*{theorem 2}{\rm\bf Proposition 2}

\theoremstyle{definition}

\theoremstyle{remark}

\def\interieur#1{\mathord{\mathop{\kern 0pt #1}\limits^\circ}}

\title[Mathematics and group theory in music]{Mathematics and group theory in music}

\date{\today}

\author{Athanase Papadopoulos}
\address{Athanase Papadopoulos, Institut de Recherche Math\'ematique Avanc\'ee, Universit{\'e} de Strasbourg and CNRS,
7 rue Ren\'e Descartes,
 67084 Strasbourg Cedex, France;  
Erwin Schr\"odinger International Institute for Mathematical Physics, 
Boltzmanngasse 9, 1090, Wien, Austria}
\email{papadop@math.unistra.fr}

\begin{document}

\maketitle


\thispagestyle{empty}


\begin{abstract}
\vskip 3mm\footnotesize{

\vskip 4.5mm
\noindent
The purpose of this paper is to show through particular examples how group theory is used in music. The examples are chosen from the theoretical work and from the compositions of Olivier Messiaen (1908-1992), one of the most influential twentieth century composers and pedagogues. Messiaen  consciously used mathematical concepts derived from symmetry and groups, in his teaching  and in his compositions.  Before dwelling on this, I will give a quick overview of the relation between mathematics and music. This will put the discussion on symmetry and group theory in music in a broader context and it will provide the reader of this handbook some background and some motivation for the subject. The relation between mathematics and music, during more than two millennia, was lively, widespread, and extremely enriching for both domains. 
\vskip 4.5mm
\noindent

\vspace*{2mm}
\noindent{ 2000 Mathematics Subject Classification: 00A65}

\vspace*{2mm}
\noindent{ Keywords and Phrases: Group theory, mathematics and music, Greek music, non-retrogradable rhythm, symmetrical permutation, mode of limited transposition, Pythagoras, Olivier Messiaen.}}

\vspace*{2mm}
\noindent{This paper will appear in the \emph{Handbook of Group actions}, vol. II (ed. L. Ji, A. Papadopoulos and S.-T. Yau), Higher Eucation Press and International Press.}

\vspace*{2mm}
\noindent{The author acknowledges support from the Erwin Schr\"odinger International Institute for Mathematical Phy\-sics (Vienna). The work was also funded by  GREAM (Groupe de Recherches Exp\'erimentales sur l'Acte Musical ; Labex de l'Universit\'e de Strasbourg), 5, all\'ee du G\'en\'eral Rouvillois
CS 50008
67083 Strasbourg Cedex. Email:  papadop@math.unistra.fr}

\end{abstract}

 \section{introduction}

\hfill  \emph{Mathematics is the sister as well as the servant of the arts.}

\hfill (Marston Morse in  \cite{Morse1959})

\bigskip

 \section{introduction}

Music is  a privileged ground for an alliance
  between arts and
sciences, and in this alliance, mathematics plays a central  role. In the first part of this paper, I will  highlight some elements of this relation and I will also point out some important works done in this area which are due to mathematicians and which are spread over several centuries. In the second part (\S 3 to 6), I will discuss in some detail the group theory that is involved in the compositions and in the theoretical work of Olivier Messiaen, one of the major twentieth-century composers and music teachers.  

Let me mention right away that besides group theory, there are many other fields of mathematics that are involved in music theory, in composition and in musical analysis: geometry, probability, category theory, combinatorics, graph theory, etc., but I do not develop any of these ideas here because they do not really belong to the subject of this handbook. In fact, any mathematical theory or idea may have its counterpart in music. Let me also mention that there are presently some very active music research groups in which mathematics plays a major role, like the IRCAM group in Paris,\footnote{Institut de Recherche et Coordination Acoustique/Musique.} whose members include Moreno Andreatta, Emmanuel Amiot, G\'erard Assayag, Chantal Buteau, Marc Chemillier, Jan Haluska, Franck Jedrzejewski and Fran\c cois Nicolas. Andreatta is one of the leaders of the group; several of his works are related to group theory; see his habilitation document \cite{Andreatta1} and the references there; e.g. \cite{A-ouvert}  \cite{AMNL} \cite{A-A} \cite{ANAAA}. The writings of Jedrzejewski are also related to symmetry and groups; see \cite{J1} \cite{J2} \cite{J3}   \cite{Jedr} \cite{J4} \cite{J5} \cite{J6} \cite{J7} \cite{J8}. There are many other modern writings on mathematics and music, done by researchers in France and outside France; see e.g.  \cite{Assayag}, \cite{Chew},
 \cite{Fau},  \cite{Fripertinger1}, \cite{Fripertinger2}, \cite{GM},  \cite{Hal},   \cite{Jo-J} \cite{Klou}, \cite {Lewin},   \cite{Tim}; see also the short introduction \cite{Papa-Intelligencer} which is more intended for mathematicians. 

I would like to thank Moreno Andreatta, Mattia Bergomi, Pierre Jehel, Thomas Flore, and Franck Jedrzejewski for their detailed comments on an early version of this paper. Andreatta and Jedrzejewski gave me (among other things) valuable biographical references. 

 \section{A brief overview of the interaction between \\
 mathematics and music}\label{s:brief}

Historically, mathematics and music are intricately linked.
Pythagoras, who is considered
as the founder of the first school of mathematics as a purely 
 deductive
science, is also the founder of a school of theoretical music (may be also the first one).\footnote{It is fair to add right away that
 the culture of Greek antiquity is, in its turn,  indebted to other cultures. Several of the major Greek philosophers and scientists
  travelled widely, and they
 acquired  an important part of their 
 knowledge from older Eastern civilizations. 
 For instance, the compiler Plutarch (c. 40-120 {\sc a.d.}) writes in his essay {\it On Isis
 and Osiris} (\cite{Plutarch-Isis} Chap.  10, 354B p. 130): ``Solon, Thales, Plato, 
 Eudoxus and Pythagoras,
  and some say,
 Lycurgus [...] came to Egypt and were in touch with the priests
 there". Plutarch even gives the names of the Egyptian priests from whom these scholars received their teaching. In the same book, Plutarch (who was a Delphic priest), provides an explanation of some of Pythagoras' aphorisms  by making a comparison between these sayings  and  Egyptian
(hieroglyphic)  writings. There are many other sources of information on the influence of Eastern civilizations on Greek culture.} 
Besides being a mathematician, Pythagoras
 was a music theorist and a composer, and his biographers describe him as 
playing  several instruments (see for instance \cite{Porphyre2} and \cite{Iamblichus}). We owe him the discovery of
the fundamental correspondence between musical intervals (that is, pairs of pitches, or of musical notes) and numerical ratios. The quickest way to describe this correspondence is by saying that to a musical interval, we associate the ratio of the frequency of the higher note to that of the lower-pitched note. Although the Pythagoreans did not talk about the frequency of a note, they were aware of this correspondence between musical intervals and fractions. The Greek mathematicians were aware of the fact that sound is produced from a periodic vibration in the air, and that a sharper note corresponds to a more rapid vibration. These views were known for instance to the  author of the \emph{Division of the canon} \cite{Euclid-d}, which is presumably Euclid, based on an earlier version due the Pythagorean mathematician Archytas (428-347 {\sc b.c.}).

We also owe to Pythagoras the first classification of consonant intervals. We recall that consonance results from playing together two (or more) different sounds, and the main question in this domain is when does such a combination give a harmonious (or ``consonant") sound and what is the reason for that. This question occupied several mathematicians and scientists, and among those who wrote on this subject we mention Aristotle, Euclid, Ptolemy, Descartes, Huygens, Galileo, Kepler, Mersenne, d'Alembert and Euler.

The two major Pythagorean discoveries, namely, the correspondence between musical intervals and numerical ratios, and the classification of consonances together with the questions related to this classification, are at the basis of all the subsequent theories of harmony. 

Pythagoras did not leave anything written -- or at least, no writing of him survives.\footnote{According to Pythagoras' biographers, it was part of his strict rules -- which he applied to himself and to his followers -- that the discoveries and the results in mathematics and music theory obtained by the members of his school should not be written up but only taught orally to the small circle of devotees which constituted that community; see e.g. \cite{Iamblichus}.} But several treatises on harmony, written by later mathematicians and based on Pythagorean ideas, survive at least in part; we shall mention a few of them below.  The discoveries of Pythagoras are described in the \emph{Handbook of Harmonics} of  Nicomachus\footnote{Nicomachus of Gerasa (c. 60-120 A.D) was a neo-Pythagorean mathematician, well known for his \emph{Introduction to Arithmetic} and his \emph{Enchiridion} (or \emph{Handbook of Harmonics}). His reverence for number is expressed in another work called \emph{Theologumena Arithmetica} (Theology of number).} (see \cite{Barker}, \cite{Nicomaque} \cite{Delatte}) and in a biography written by  Iamblichus \cite{Iamblichus}.\footnote{Iamblichus (c. 245-325 {\sc a.d.}) was a neo-Platonist philosopher, known for his cosmologica system based on mathematical formalism. Both Nicomachus and Iamblichus were Syrian. (Syria,  at that time, was a Roman province). Only a small portion of the works of Nichomachus and of Iamblichus survives.} Aristotle, who is a reliable source, reports in  his \emph{Metaphysics} (\cite{Aristotle-Metaphysics} A5, 986a16) that Pythagoras used to say that
   ``everything  is number".

The works of the Pythagoreans reached us in the form of quotations, in relatively small number, but very rich in content, see e.g. the volumes \cite{Presocratiques}, \cite{Barker}, and \cite{Pres}.

Let us also recall that, in principle, every mathematical treatise of classical Greece 
contained  a chapter on music. In fact, such a treatise usually consisted of four parts: Number theory, Music, Geometry, and Astronomy, in that order, because the part on music was based on the results of number theory, and the part on astronomy was based on the results of geometry. To give the reader an idea of the important connections between number theory and music, we recall that the theory of proportions and the theory of means were developed precisely for their use in music theory. The division of a musical interval into two or more subintervals is defined in terms of proportions that depend on the pitches of the musical notes involved, and this division was formulated in terms of ratios of lengths of the subintervals to the length of the whole interval. This mathematical theory of division of musical intervals was made possible by the Pythagorean discovery of the correspondence between musical intervals and fractions that we mentioned, and of the logarithmic law that governs this correspondence (concatenation of musical intervals corresponds to a multiplication at the level of the numerical values).\footnote{\label{f:6}This logarithmic law was known in Greek antiquity, and it was used long before logarithms were formalized by mathematicians. 
  For instance, Theon of Smyrna in his \emph{Mathematical exposition} writes in \cite{Theon} (p. 103): ``Since the  ratio of the consonance of octave is $2/1$ and the one of the consonance of  a fourth is $4/3$, the ratio of their sum is $8/3$." That is, he knew that the ratio of the sum of two musical intervals is the product of the corresponding ratios. Likewise, Iamblichus, in \cite{Iamblichus} \S 115, reports that Pythagoras noticed that ``that by what the fifth surpasses a fourth is precisely the ratio of $9/8$"; in other words, he saw that the difference between the two ratios $3/2$ and $4/3$ is $9/8$ (which corresponds to the fact that $9/8$ is $3/2$ divided by $4/3$).}

To introduce things more precisely, let us give a concrete example of how music theory acted as a motivation for number-theoretic research. The example concerns the arithmetic of the so-called superparticular ratios. These are the numerical ratios of the form $(n+1)/n$, where $n$ is a positive integer. Superparticular ratios are important for several reasons, one of them being that the corresponding musical intervals appear as the successive intervals in the decomposition of a sound into harmonics, and thus we constantly hear them in any sound that is produced around us. Therefore, the ear is familiar with them, and this makes them important. Another reason for which these ratios are consequential is that in the ancient Greek classification, the so-called \emph{consonant} intervals are either superparticular or of the form $n/1$ (see e.g. \cite{Theon}, Chapter II). We already mentioned (see Footnote \ref{f:6}) the values of the octave ($2/1$), the fifth ($3/2$), the fourth ($4/3$) and the major tone ($9/8$). Some other important superparticular ratios which are useful in music appear in the tables that follow.  The major and minor thirds, defined respectively by ($5/4$) and ($6/5$), started to be considered as ``imperfect consonances" at the thirteenth century and they played, after that period an important role in composition. The ``Didymus comma" (also called the syntonic comma), whose value is  ($81/80$), is the difference between the major tone ($9/8$) and the minor tone ($10/9$), and it played a significant role in Greek theoretical music.

Besides the question of classification, which was motivated by music, there are purely mathematical developements. Indeed, several natural questions concerning superparticular ratios were formulated and studied by the mathematically-oriented Pythagoreans. I shall mention a few of them as examples; some of these questions are easy, and others are difficult.

   \begin{enumerate}
 \item \label{q1} Can the square root of a superparticular ratio be superparticular?  Or at least, can it be rational? 
 
 A musical naturally related question is the following: can we divide a consonant interval into two equal consonant intervals?

  \item \label{q2} Given a superparticular ratio, can we enumerate all the various ways of expressing it as a product of superparticular ratios?  Is this number finite or infinite? Is the number of possibilities finite if we fix a bound on the number of factors?
 
 This question is related to the question of dividing a consonant interval into a certain number of consonant intervals and the problem of constructing scales whose all intervals are consonant.

 \item \label{q3} Given a finite set of primes, e.g. $\{2,3\}$, $\{2,3,5\}$ or $\{2,3,5, 7\}$, is the set of superparticular ratios whose prime factors (of the numerator as well as the denominator) belong to that set always finite? Can we enumerate all the elements of this set? 
 
 This question is related to the construction of scales out of a finite set of prime numbers. As an example, the reader can notice that the numerical values 
 of the ratios that appear in the following table extracted from Descartes'  \emph{Compendium} are all multiples of 2, 3 and 5.\footnote{Descartes writes in \cite{Descartes-Compendium}, p. 122: ``All the variety of sounds, for what concerns pitch, 
 stems only from the numbers 2, 3 and 5; and all the numbers that define the [musical] degrees as well as dissonance are multiples of these three sole numbers".} 
 \end{enumerate}

 \begin{table}[h]
 \begin{center}
 \renewcommand{\arraystretch}{1.25}
 \begin{tabular}{ | l | l | l | l | l | l | l | l | l | l | l  }
 \hline
$2/1$ & 8ve  &&&& &&&& \\ \hline
$3/1$ & 12th  & $3/2$  &  5th &&&  &&&   \\ \hline
$4/1$ & 15th &$4/2$   & 8ve & $4/3$ &   4th &&  &&\\ \hline
$5/1$ & 17th & $5/2$  & 10th maj. & $5/3$ & 6th maj. & $5/4$ & 2nd &   &\\ \hline
$6/1$  & 19th & $6/2$ & 12th & $6/3$ & 8ve& $6/4$ & 5th & $6/5$ &3rd min.  \\ \hline
 \end{tabular}
 \medskip
 \caption{A table of musical intervals, ordered according to the denominator; extracted from Descartes' {\it Compendium musicae}, \cite{Descartes} Tome X p. 98.}
  \label{Table:Descartes1}
 \end{center}
\end{table}

 Let me make a few more comments on these questions.

 The response to Question \ref{q1} is known since antiquity. A proof of the fact that there is no rational fraction whose square is equal to a superparticular ratio, attributed to Archytas, a Pythagorean from the first half of the fourth century {\sc b.c.},  is contained in Boethius' \emph{Musical Institution}, Book III. A more general result is a consequence of Proposition 3 of the Euclidean \emph{Section of the Canon} which says the following: \emph{For any pair of integers  $B,C$ whose quotient is equal to a superparticular ratio,  there is no sequence of integers $D, E,F,\ldots N$ between $B$ and $C$ satisfying $B/D=D/E=E/F=\ldots=N/C$.}
   (cf. \cite{Barker}, vol. II, p. 195). Note that in this statement, the fraction $B/C$ is not necessarily in reduced form.

 In Boethius' \emph{Musical Institution} III.5 d \cite{Boethius}, the author mentions that in order to circumvent the impossibility of dividing the tone $(9/8$) into two equal parts, Philolaos divided it into two unequal parts, the one being ``less than a semitone", which he called the \emph{diesis} or \emph{lemma}, and which is also called the \emph{minor semitone}, and  the other one being ``greater than a semitone", called \emph{apotome}. The \emph{comma} is the difference between these two intervals. (See also \cite{Presocratiques} p. 500).

  Regarding Question 2, one can note that any superparticular ratio can be written as a product of two others, using the following: 
  \[\frac{p+1}{p}=\frac{2p+2}{2p}= \frac{2p+2}{2p+1}\times \frac{2p+1}{2p}.\]
We can then apply the same trick to the fraction $\displaystyle \frac{2p+2}{2p+1}$, or to $\displaystyle \frac{2p+1}{2p}$ (or to both), and therefore the process goes on indefinitely. In particular, this shows that every superparticular ratio can be written as a product of superparticular ratios in an infinite number of ways.

 Aristides Quintilianus, in Book III of his {\it De Musica} \cite{Aristides}, used this method to describe the division of the tone as $17/16\times 18/17= 9/8$, and then the following division of semitones and of quarter tones:   
$33/32 \times 34/33= 17/16$ and $35/34\times 36/35=18/17$. 

There is another general method for obtaining products of superparticular ratios which is based on the following equalities:
\[\frac{p+1}{p}=\frac{3p+3}{3p}= \frac{3p+3}{3p+2}\times \frac{3p+2}{3p+1}\times \frac{3p+1}{3p}.\]
This gives the following well known division of the fourth:
 \[\frac{4}{3}= \frac{12}{11}\times \frac{11}{10}\times \frac{10}{9}.\]

Question 3 was solved in the affirmative by the mathematician Carl St\o rmer who proved in \cite{Stormer} (1897) that for any finite set of primes $\{p_1,\ldots,p_n\}$, there are only finitely many superparticular ratios whose numerator and denominator are products of elements in this finite set. He also described a procedure to find such fractions. As a consequence, Table 3 gives the list of all superparticular ratios whose prime factors of the numerator and denominator belong to the set $\{2, 3, 5\}$. Let us note that all the fractions in Table 3 were used in music, since Antiquity. This is another example where music theorists were far ahead of mathematicians. For more on superparticular ratios in music, see \cite{Halsey-Hewitt}. See also \cite{Read} and \cite{Reiner} for an account of some combinatorial problems related to music theory.

Our next example is extracted from Ptolemy's \emph{Harmonics} (see \cite{Barker} p. 203 and \cite{Presocratiques} p. 528), where the author comments on the following divisions of the fourth  ($4/3$) (such divisions are traditionally called tetrachords), which are due to Archytas, and which are called respectively \emph{enharmonic}, \emph{chromatic} and \emph{diatonic}:
 $${5\over 4}\times
{36\over 35}\times{28\over 27}=
 {32\over 27}\times {243\over 224}\times{28\over 27}=
 {9\over 8}\times {8\over 7}\times{28\over 27}={4\over 3}.$$

Several tables of ancient Greek musical 
   tetrachords are contained in  Reinach's  {\it Musique Grecque},  \cite{Reinach}  see e.g. Table \ref{Rei}. Most (but not all) of the numerical values in these tables are superparticular ratios.
 
 \begin{table}[h]
 \begin{center}
 \renewcommand{\arraystretch}{1.25}
 \begin{tabular}{ | l | r | r | r  | r | r |}
 \hline
 Archytas  &  9/8 & 8/7 & 28/27 \\ \hline
 Eratosthenes  & 9/8   & 9/8 & 256/243 \\ \hline
 Didymus  &  9/8 & 10/9 & 16/15 \\ \hline
 Ptolemy  &  10/9 & 11/10 & 12/11 \\ \hline
 \end{tabular}
 \medskip
 \caption{Diatonic genus (after Reinach)}
  \label{Rei}
 \end{center}
\end{table}

Some more questions on superparticular ratios in music are discussed in \cite{Halsey-Hewitt}

The next table of intervals is extracted from Euler's book on music, the \emph{Tentamen} \cite{Eul01}, which we shall mention again in what follows. In the tradition of the Greek musicologists, Euler made a systematic classification of the useful musical intervals according to their numerical values and he developed a theory of the musical significance of the ordering in these lists. His tables, in the  \emph{Tentamen}, involve the prime numbers 2, 3 and 5. But Euler also used the set $\{2, 3, 5,7\}$, for instance in his memoir {\it  Conjecture sur la raison de quelques dissonances g\'en\'e\-ra\-lement re\c cues dans la musique} \cite{Eul02},  and this was considered as a novelty, compared to the smaller set $\{2, 3, 5\}$ which was used by his predecessors in the post-Renaissance Western world.

In the twentieth century,
Hindemith, in his famous treatise \cite{Hindemith}, 	also uses the integer 7. In Greek antiquity not only the number 7 was used, but in principle no number was excluded; we already mentioned a few examples. Aristoxenus  of Tarentum, the great Greek music theorist of the fourth century {\sc b.c.}, tried to make exhaustive lists of scales where a large number of primes appear, see \cite{Belis}. 
 
 \begin{table}[h]
 \begin{center}
 \renewcommand{\arraystretch}{1.25}
 \begin{tabular}{ | l | l |}
 \hline
$2/1$ &  Octave   \\ \hline
$3/2$  & Fifth   \\ \hline
$4/3$ & Fourth  \\ \hline
$5/4$ &  Major Third   \\ \hline
$6/5$ &  Minor Third   \\ \hline
$9/8$ & Major Tone   \\ \hline
$10/9$ &  Minor Tone  \\ \hline
$16/15$ &  Diatonic Semitone   \\ \hline
$25/24$  & Chromatic Semitone   \\ \hline
$81/80$ &  Didymus Comma   \\ \hline
 \end{tabular}
 \label{Table:Stormer}
 \medskip
 \caption{The list of superparticular ratios whose prime factors belong to the set $\{2, 3, 5\}$.}
 \end{center}
\end{table}
  
Some relations between music theory and number theory are also  manifested by the terminology. The Greek  word ``diastema" means  at the same time ``ratio" and ``interval". The same is true for the word ``logos".\footnote{Theon's treatise \cite{Theon} contains a section on the various meanings of the word ``logos".} The theory of means, in Ancient Greece, found its main applications in music. Defining the various means between two given integers $a$ and $b$ ($a<b$) was seen practically as  inserting various notes in  the musical interval whose numerical value is the quotient $b/a$. For instance, the harmonic mean  of the interval $[6,12]$ (which is 8) corresponds to the note that divides an octave into a fourth followed by a fifth. Thus, it is not surprising that the oldest expositions of the theories of proportions and of means are contained in musical textbooks, and the examples, in these writings, that illustrate these mathematical theories are often borrowed from music theory.\footnote{Examples of computations illustrating mathematical theories that have a musical significance may also be found in later works. For instance, in his famous \emph{Introductio in analysin infinitorum}  (Introduction to the analysis of the infinite, published in 1748),  Euler, while presenting  his methods of computation using logarithmes, explains how one can find the twelfth  root of 2, which in fact is the value of the unit in the chromatic tempered scale. In Chapter VI of the same treatise, Euler works out an approximate value of $2^{7/12}$, which of course corresponds to the fifth. There are other examples of this sort.} The discovery of irrational numbers was motivated in part by the mathematical difficulty of dividing a tone into two equal parts. The distinction continuous vs. discontinuous arose from the attempt of splitting up the musical continuum into the smallest audible intervals. We note by the way that not all the intervals useful in music were rational. Aristoxenus made a distinction between rational and irrational musical intervals. 

 There are also important repercussions of musical theories in geometrical problems, e.g. on the geometric divisions of the musical intervals and on the geometric constructions on means. Ptolemy (c. 90-168 {\sc a.d.}), in his \emph{Harmonics} (Book II, ch. 2) describes a geometric instrument, called \emph{helicon}, which was used to measure consonances. 
There were also impacts on famous problems like the duplication of the cube and on several questions on constructions with compass and straightedge. This came very naturally, since the same people who worked on music theory worked on these geometrical problems.

The division of the teaching of mathematics into four parts, which was given later on the name {\it Quadrivium}\footnote{The Latin word \emph{quadrivium} was introduced by Boethius (5th century {\sc a.d.}).}  (the ``four ways") lasted  until the middle ages, and the status of
theoretical music as  part of mathematics persisted in Western Europe  until
the beginning of the 
Renaissance (c. 1550). A textbook on the quadrivium available in French translation \cite{Theon} is the one written by Theon of Smyrna (c70-135 {\sc a.d.}) which we already mentioned.\footnote{The book, in the form it survives, contains three parts; the part on geometry is missing.}
 
One of the oldest Pythagorean texts that survives  describes  geometry, arithmetic, astronomy (referred to as \emph{spherics}), and music 
  as ``sister sciences". This text is a fragment from a book titled {\it On mathematics } by Archytas, a Pythagorean from the first half of the fourth century {\sc b.c.} and it is known through a quotation by the   philosopher   Porphyry\footnote{Porphyry (c. 233-309 {\sc a.d.}) was a Hellenized Phoenician, born in Tyre (presently in Lebanon). In 262 he went to Rome, where he stayed six years, during which he studied under Plotinus, one of the main founders of neo-Platonism. He is known for his {\it Commentary on Ptolemy's Harmonics}  and for several books on philosophy and a book on the history of philosophy. His \emph{Pythagorean life} is part of the latter. He wrote a \emph{Life} of his master Plotinus, and he  edited his works under the name of \emph{Enneads}.} in his {\it Commentary on Ptolemy's Harmonics} (part of which is attributed to the mathematician Pappus) \cite{Porphyre1}. The text of Porphyry was later on 
   edited, with a Latin translation accompanying the Greek original, by the mathematician John
   Wallis (1616-1703), and it was published as part of his collected works  ({\it Opera
  Mathematica} \cite{Wallis-Opera} Vol. III).\footnote{Wallis also worked on critical editions of Ptolemy's
    {\it Harmonics}  and of the
    {\it Harmonics}
  of Manuel
    Bryennius, the fourteenth century Byzantine music theorist.  
     These two 
    editions, together with teh one of Porphyry's {\it Commentary to   Ptolemy's   Harmonics}, with 
    a Latin translation accompanying the Greek text, together with editions of
     works 
     by Archimedes and Aristarchus of Samos,
      constitute Volume III of Wallis' collected works, published  
      in  three volumes  
    in Oxford, in  1699.} Porphyry writes: 
  \begin{quote}\small
  Let us now cite the words of Archytas the
  Pythagorean, whose writings are said to be mainly authentic. In his book {\it
  On Mathematics}, right at the beginning of the argument, he writes:  ``The
  mathematicians seem to me to have arrived at true knowledge, and it is not
  surprising that they rightly conceive the nature of each individual thing; 
  for, having reached true knowledge about the nature of the universe as a
  whole, they were bound to see in its true light the nature of the parts as
  well. Thus, they have handed down to us clear knowledge about the speed of the
  stars, and their risings and settings, and about geometry, arithmetic, and
 spherics, and, not least, about music; for these studies appear to be
  sisters".\footnote{This English translation is taken from  the   {\it Selections illustrating the history of Greek Mathematics}, edited by
    Ivor Thomas, see \cite{Thomas}, Vol. 1, p. 5. The text is also quoted in French in the volume \cite{Presocratiques}, p. 533.}
    \end{quote}
    
 The use of the word ``sister" in the preceding quote is similar to the one in the quote by Marston Morse (whom we shall mention again below) which is at the beginning of the present paper.

Euclid  wrote
 several  treatises on music. Among them is the
 {\it Division of the canon} which we already mentioned, in which he gives an account of the Pythagorean
 theory of music, and  
 which contains in particular  a careful exposition of the mathematical theory of proportions applied to musical harmony.
 Proclus\footnote{Proclus (412-485 {\sc a.d.}) first studied mathematics in
    Alexandria under Heron, and then  
    philosophy in Athens under Plutarch. He became the head of the neo-Platonic school of Athens, after Plutarch and Syrianus.} in his 
 {\it Commentary to Euclid's  First Book of Elements} attributes to Euclid another treatise titled 
 {\it Elements of Music}, which unfortunately
 did not survive into our time. 
 Eratosthenes (c. 276-194 {\sc b.c.}) also had   
  important impacts  on both  fields, mathematics and music.
  His work, the {\it Platonicus},
  contains a section on music theory  which is referred to several times
  by Theon of Smyrna in \cite{Theon}.

  Several of the mathematicians-musicians we mentioned
 were equally erudite in other domains of knowledge. For instance, Eratosthenes, who 
 was the administrator of the famous library of Alexandria, was considered as the most learned person of his time, and for this reason he
  was known under the name ``$\beta$", the second letter of the alphabet, which was a manner  of indicating
  that he was ``second" in every domain of knowledge. Ptolemy, whom we already mentioned, was a
 mathematician, geographer, astronomer, poet and expert in oriental mysticism, and he was probably the greatest music theorist of the Greco-Roman period. His major work, the  {\it
Mathematiki Syntaxis} ({\it Mathematical collection}), a treatise on astronomy in 13 books,  reached us through the Arabs with the title {\it
Almagest} (a corrupt form of the Greek superlative Megistos, meaning ``the greatest").  Ptolemy is also the author of an important   musical treatise, 
the {\it Harmonics}, in which he exposes and 
 develops    Pythagorean musical theories.
  This treatise was also translated
 into  Arabic
   in the ninth century and into Latin in the sixteenth century.
   From the later Greek period,
  we can mention the
   mathematician Pappus  
   (third   century  {\sc a.d.})
  who, like Euclid, Eratosthenes, and Ptolemy, lived in
   Egypt, and who also  was   an excellent  music theorist. He wrote
   an impressive exposition of 
    all of what was known in
     geometry at the time, 
     with the title  {\it Synagoge} (or {\it Collection}). He  
    also wrote
   a Commentary on Euclid's {\it Elements} and a  
      Commentary on Ptolemy's work, including his {\it Harmonics}.
   Proclus, whom we already mentioned several times and the author of
   the famous {\it Commentary to Euclid's First Book of Elements},  also wrote commentaries on several of Plato's dialogues, including the
   {\it Timaeus}, a dialogue which is essentially a treatise on mathematics and music. The subject of this dialogue is the creation of the universe, described allegorically as a long musical scale. 

    The
belief in a strong connection between the four fields of the quadrivium is also part of a broader deep feeling of order and harmony in nature and in human kind.   This is also the origin of the word ``Harmonics", which is used in many places instead of the word music. This word
  has the flavour of  order, of structure  and of measure. This feeling of harmony which was
    shared by most of the major  thinkers   of Greek antiquity  was  a vehicle for an extraordinary flourishing of arts and sciences which included  the development of an abstract and   high-level
 mathematics  and the construction of coherent
 systems governing  the sciences of  
  music, astronomy, physics,  metaphysics, history and theatre. 
  It is generally accepted  that these systems 
  had a real  and probably irreversible  impact 
  on all
   human thought and in any event, they continued to dominate most branches of knowledge in  Europe 
  until the end of the middle ages.

The belief in an intimate  relation between mathematics and music, which was stressed primarily by the Greek thinkers, sometimes took the
form of a belief in the fact that music -- not only its theory, but also the emotion that it produces -- is in many ways  identical to the emotion that mathematical pure thought can produce. Such a feeling was also formulated in modern thought. Let us quote for instance Marston Morse, from his paper  \emph{Mathematics and the
Arts} \cite{Morse1959}: 
\begin{quote}\small
Most convincing to me of the spiritual relations between mathematics
and music, is my own very personal experience. Composing a little in an
amateurish way, I get exactly the same elevation from a prelude that has come to
me at the piano, as I do from a new idea that has come to me in mathematics.
\end{quote}

Although, by the end of  the sixteenth century,  
the antique tradition considering music as 
part of mathematics progressively disappeared,
the development of music theory and practice continued to be accompanied with a fruitful alliance  with mathematics. 

Among the seventeenth century mathematicians involved in this alliance, we first mention Newton, one of the principal founders of modern science.

 Newton was interested in every kind of intellectual activity, and of course he was naturally led to music theory.\footnote{Pythagoras, for whom Newton had a great respect, is mentioned several times in the \emph{Principia}. In the Scholia on Prop. VIII Book III on universal gravitation, Newton declares that Pythagoras was aware of several physical laws, for instance the fact that square of the distance of the planets to the sun is inversely proportional to the weights of these planets, but that because of the nature of his teaching (which was essentially esoterical), nothing written by him could survive.}
 A notebook left from his early college days (c. 1665) concerns this subject, and it contains, in the old Greek tradition, a theory and computations of the division of musical intervals. Newton is also known for the use of the logarithms in his musical computations. He discussed several points of music theory in his correspondence, in particular with John Collins \cite{Newton-Corresp}. He made relations between some
divisions of musical intervals, and in particular of symmetrical divisions (palindromes), and questions in optics on the division of
the color spectrum. In 1666, Newton discovered that sunlight is a mixture of several colors, and this was one of the starting points for his theory on the correspondence between the color spectrum and the musical scale, which became later  one of his favorite subjects. This topic is also discussed in his correspondence \cite{Newton-Corresp}, in particular with Henry Oldenburg in 1675 and with William Briggs in 1685 and in his popular work  \emph{Opticks} (1704). Newton's theory of sound is also discussed in his famous paper \emph{New Theory about Light and Colors.} (1672). Let us note right away that the correspondence between colors and pitch is also one of the main themes in the theoretical work of Olivier Messiaen that we shall discuss later in this paper.

 One should emphasize here that Newton's ideas about the relation between the spectrum of colors and the  musical diatonic scale, and more precisely the fact that the two spectra are governed by the same numerical ratios, was part of his firm convection that the same universal laws rule all aspects of nature.  
 Voltaire was one of several theorists on the continent who were eager to adopt and promote Newton's ideas, in particular his theory concerning the relation between the seven-scale color spectrum and the seven-scale diatonic scale, see \cite{Voltaire-Elemens} and \cite{Voltaire-Table}.
 
  After Newton, it is natural to mention Leibniz, with his famous sentence: \emph{Musica est exercitium arithmeticae occultum nescientis se numerare animi} (Music is a secret arithmetic exercise of the mind which is unaware of this count), that is, music consists in a mathematics count, even though his listener is unconscious of that. The sentence is extracted from the correspondence of Leibniz with Christian Godlbach   \cite{LB1} \cite{LB}, the famous number theorist and friend of Euler.

Without going into any detail, we now mention  some seventeenth century works on music written by scientists. 
Kepler's famous \emph{Harmony of the world} (1616) \cite{Kepler} contains several sections on music theory, written in the Pythagorean tradition. (We note by the way that Kepler described himself as a ``neo-Pythagorean"). See also \cite{Papa-Kepler} on the relation between mathematics and music in Kepler's \emph{Harmony of the world}. The first book that Descartes wrote is a book on music, \emph{Compendium Musicae} (1618) \cite{Descartes-Compendium}. Mersenne, the well known number theorist, wrote a music treatise called \emph{Trait\'e de l'harmonie universelle} (1627) \cite{Mersenne}. In this treatise, he states, on p. 35: ``Music is part of mathematics", on 39: ``Music is a science; it has its real proofs which are based on its proper principles", and on p. 47: ``The music I consider is subordinate to arithmetic, geometry, and physics". Galileo Galilei's \emph{Discourses and dialogues concerning the two sciences} (1638) \cite{Galileo} , which was his last writing, contains sections on theoretical music.
 Christiaan Huygens also wrote important treatises on music, e.g. his \emph{Letter concerning the harmonic cycle} \cite{Huygens-cycle} and his works on multi-divisions of the octave \cite{Huygens1}. There are several other works of Huygens on music theory. There are also several sets of letters on music theory in the correspondence of several mathematicians, including Descartes, Huygens and Leibniz,\footnote{We recall that in that period there were still very few scientific journals, and that scientists used to communicate their results by correspondence. The letters of major mathematicians were collected and published, usually after their author's death, but sometimes even during their lifetime. The correspondence \cite{Descartes-C} of Descartes occupies Volume 1 to 5 of his twelve-set Collected Works \cite{Descartes}. The correspondence of Euler occupies several volumes of his Collected Works \cite{Euler}, and up to now only part of it has been published.} and we already mentioned the correspondence of Newton \cite{Newton-Corresp}. 

Among the eighteenth century mathematicians who worked on music, we mention 
Euler, who wrote a book, \emph{Tentamen novae theoriae musicae ex certissimis harmoniae principiis dilucide expositae} (Essay on a new musical theory exposed in all clearness 
according to the most well-founded principles of harmony), already mentioned,  and several memoirs on music theory,\footnote{There are also several papers of Euler on acoustics, but this is another subject.} see \cite{Eul01}, \cite{Eul02}, \cite{Eul033} \cite{Eul04}. See also the forthcoming books \cite{HP2014} and \cite{Euler-Hermann} on Euler's musical works. Euler formulates as follows the basic principle on which he builds his music theories:
\begin{quote}\small
 What makes music pleasant to our ears depends neither on  will nor on habits. [...] Aristoxenus denied the fact that one has to search for the pleasant effect of music in the proportions established by Pythagoras ... Led by reasoning and by experiments, we have solved this problem and we have established that two or more sounds produce a pleasant effect when the ear recognizes the ratio which exists among the number of vibrations made in the same period of time; that on the contrary their effect is unpleasant when the ear does not recognize this ratio (Extracted from the Introduction of \cite{Eul01}; see also the French translation in \cite{Euler-Hermann}).
 \end{quote} 
 
In fact, the book \cite{Eul01} is the first one which Euler wrote. He finished writing it in 1731, the year he obtained his first position, at the Saint-Petersburg Academy of Sciences, and he was 24 years old. It is most probable that several projects in combinatorics and in number theory occurred to Euler while he was developing his music theory, since several natural questions  regarding primes and prime factorizations of numbers appear in that theory.  
  
 It is also interesting to hear what composers say about mathematics. 

We can  quote Jean-Philippe Rameau, the great eighteenth's century French
composer and music theorist, from his famous  {\it
Trait\'e de l'harmonie r\'eduite \`a ses principes naturels} (1722) (see
\cite{Rameau-Traite}, Vol. 1, p. 3)\footnote{I am translating from the French.}:
\begin{quote} \small
Music  is a science
which must have determined rules. These rules must be drawn from a
principle which should be evident, and this principle 
cannot be known 
without the help of mathematics.
I must confess 
that in spite of all the experience which I   have acquired in music by
practicing  it for a fairly long period of time, 
it is nevertheless only 
with the help of mathematics that my ideas became
  disentangled
and that light  
succeeded
to a certain darkness of which I was not aware before.
\end{quote}

   In his {\it
D\'emonstration du principe de l'harmonie},
Rameau relates  how, since his childhood, he was aware   of
the role that mathematics  plays in music (\cite{Rameau-Traite}, Vol. 3, p. 221):
\begin{quote} \small
Led, since my early youth, by a mathematical instinct in the study of an art
for which I found myself destined, and which occupied me all my life long, I wanted to know its true principle, as the only way to guide me with certitude, regardless of the problems and accepted ideas.
\end{quote}

Two hundred years after Rameau,
Olivier Messiaen
made similar statements concerning the relation between mathematics and music, and we shall record them in \S \ref{s:OM} below.

Rameau wrote a major corpus of works on music theory. They include  his \emph{Trait\'e de l'harmonie r\'eduite \`a ses principes naturels} \cite{Rameau-Traite}, his \emph{Nouveau Syst\`eme de Musique Th\'eorique} \cite{Rameau-nouveau}, his  \emph{D\'emonstration du principe de l'harmonie} \cite{Rameau-Demonstration} and there are many others; see the whole collection in \cite{Jacobi}. There is also a correspondence between Euler and Rameau, see \cite{Euler-Hermann} and \cite{HP2014}.

In a review of Rameau's \emph{Trait\'e de l'harmonie r\'eduite \`a ses principes naturels} which he wrote in the famous \emph{Journal de Tr\'evoux}, the Jesuit mathematician and philosopher L.-B. Castel wrote: ``Music is henceforth a vast quarry which will not be exhausted before a long time, and it is desirable that philosophers and geometers will want to lend themselves to the advancement of a science which is so puzzling."

D'Alembert also became very much interested in Rameau's theoretical writings, and he wrote an essay explaining his theories \cite{Alembert}. The relation between the two men became tense and eventually bad -- d'Alembert accused Rameau of exaggerating the role of mathematics in his music -- but this is another story.

Diderot, one of greatest figures of the French \emph{Enlightenment} and one of the two main editors (the other one being d'Alembert) of the famous \emph{Encyclop\'edie}, wrote a book on the theory of sound\footnote{In 1784, Diderot published a collection of 5 memoirs under the general title \emph{M\'emoires sur diff\'erents sujets de math\'ematiques}; three of these memoirs concern sound and music theory.} in which he writes (\cite{Diderot} p. 84): ``The musical pleasure lies in the perception of ratios of numbers [...] Pleasure, in general, lies in the perception of ratios". 
     
Finally, let us give a few examples from the modern period, by quoting a few geometers.

In a letter to his friend G. Wolff, a teacher at the conservatory of Leipzig, Beltrami writes (\cite{Beltrami-Boi} p. 154, note 111): 
\begin{quote}\small
Between music and mathematics, there is a reconciliation which has not yet been noted [...] A mathematical reasoning is comparable to a sequence of chords [...] and the discovery of a new branch of mathematics is comparable to a harmonic modulation.
\end{quote}

  In his paper \emph{Twentieth century mathematics} (1940) \cite{Morse1940}, Morse  writes the following: 
  \begin{quote}\small
  Mathematics is both an art and a science, and the lack of appreciation of this fact is responsible for much misunderstanding [...] Objective 	advances must be revised in form to make them aesthetically acceptable and logically comprehensible, while advances of a more subjective nature, if complete and harmonious, will not long remain unapplied.
  \end{quote}
   In a letter \cite{Morse1941b}, he writes:
  \begin{quote}\small
Mathematics is both an art and a science ... Mathematicians are the freest and most fiercely individualistic artists. They are subject to no limitations of materials or instruments. Their direction at any time is largely determined by their tastes and intellectual curiosity. Their studies are really the studies of the human mind. To me the work of Einstein is even more important as a free and beautiful expression of the creative imagination of an individual than as a part of the science of physics.
\end{quote}

The paper by Birkhoff at the Bologna ICM conference (1928) concerns mathematics and arts \cite{Birkhoff1929}, and there are several papers by him on the relation between mathematics and music, see e.g. his treatise \emph{Esthetic measure} \cite{Bir-E}.

The composer Milton Babbitt, whom we shall mention below and who taught music theory and sometimes mathematics at Princeton, insisted on a rigorous and a mathematically inclined teaching of music. He writes in \cite{Babbitt-congress} that a  musical theory
should be ``statable as a connected set of axioms, definitions, and theorems,
the proofs of which are derived by means of an appropriate logic". This is in the tradition of the great musical treatises  of Mersenne \cite{Mersenne}, Rameau  \cite{Rameau-Traite}, etc. which we already mentioned and where the exposition is in the form of Theorem, Lemma, etc.

Several serial composers, before Babbitt, knew that they were dealing with groups. For instance, the Austrian composer Hanns Jelinek (1901-1969), who had been in contact with Arnold Schoenberg and Alban Berg, wrote a book on twelve-tone music \cite{Jelinek}, in which he explicitly cites the "Gruppenpermutation" (p. 157, vol. 2). Likewise, Herbet Eimert, who made in 1964 a compete list of the ``all-intervall series", knew since the 1950s that he was dealing with groups. One may also cite Adriaan Fokker  (1887-1972), a Dutch physicist\footnote{Fokker obtained his PhD under Lorentz, and he also studied under Einstein, Rutherford and Bragg. His name is attached to the Fokker-Planck equation.} and musician,  who wrote extensively on music theory, cf. for instance \cite{Fokker1} and 
\cite{Fokker2}. Fokker was very much influenced by the works of Huygens on music theory.  An important article on the systematic use of group theory in music is \cite{HH} by Halsey and Hewitt.

It is always good to see, in skimming through these papers and books, how mathematics can serve
music and vice versa

\medskip

In the rest of this survey, I will concentrate on the work of the French composer Olivier Messiaen, because it involves in several aspects group theory. Messiaen stated explicitly (like Rameau, Euler and others did before him) that what makes the charm of a musical piece is the mathematical structures that stand behind it.

\section{The music of Olivier Messiaen}\label{s:OM}

We are concerned in the following pages with certain mathematical aspects of the
musical compositions and of the theoretical writings of Messiaen.\footnote{I already reviewed some of these ideas in my paper \cite{2003b} (2003).} 
Even
though Messiaen never considered  himself as a mathematician, he granted
to mathematics  a prominent place, both in his compositions
  and in his theoretical teaching.
  The titles of some of his pieces, like
 {\it Le Nombre L\'eger} (Pr\'elude No. III for piano), {\it Soixante-quatre
dur\'ees} (Piece No. VII of his {\it Livre d'Orgue}),
are significant in
this respect. 
The
mathematical notions that are involved in his compositions
  are basic notions
(permutations, symmetries, prime numbers, periodicity, etc.), and it may be worthwhile to
stress right away the fact that the fact that Messiaen 
uses these notions in a mathematically elementary and simple way does not reduce the place of mathematics in his work.  Questions related to properties of
sequences of numbers, of their transformations and of their symmetries,
however elementary they are, are part of mathematics. Messiaen worked with these notions 
consciously and systematically. In a book of dialogues  with Claude Samuel
\cite{Samuel}, he recalls that since he was a child, he was
fascinated by certain properties of numbers, properties which were led to
play a central role in his musical language. 
 In \cite{Samuel}, p. 118,
  answering a question about {\it de\c ci-t\^ala}
Indian rhythms and more generally about the reasons which guided him in his choice of certain rhythmical formulae, Messiaen says: 
\begin{quote}\small
  I was oriented
towards this kind of research, towards asymmetrical divisions, and
towards an element which one encounters in Greek meters and in Indian
rhythms: prime numbers. When I was a child, I already liked  prime
numbers, these numbers which, by the simple fact that they are not
divisible, emit an occult force [...]
\end{quote}

One  aspect of the music of Messiaen  is  a 
balance between reason and intuition, between poetic creation 
 and a rigorous formal structure. 
His theoretical work is in the tradition of the Greek quadrivium, and we can quote here 
the composer Alain Louvier, 
  who was a student of
  Messiaen at the Conservatory of Paris and who
  says, in his foreword to Messiaen's {\it Trait\'e de rythme, de couleur et
  d'ornithologie} \cite{Messiaen-Traite},  that in his teaching, Messiaen
  placed Music at the confluence of a new \emph{Septivium}:
   Mathematics, Physics,
  Cosmology, Acoustics, Physiology, Poetry, and Philosophy. Understanding the way in which 
 mathematical structures 
 are present in Messiaen's music can at least
  serve the purpose of
   making his music less   enigmatic
  than it   appears
  at first hearing. 
  
Finally,  beyond the description of
Messiaen's work, one of the themes which we would like to develop in the next sections is
that music (and in particular rhythm) is   a certain way of
giving life to mathematical structures, and of rendering them perceptible to our senses. More than that, music  transforms these notions into
   emotionally
affecting objects.

We have divided the rest of our exposition into three parts, with the following titles:
 \begin{itemize}
 
\item Rhythm.

\item Counterpoint.

\item Modes of limited transposition.

\end{itemize}
In each part, the reader will notice the relation with group theory.


\section{Rhythm}\label{s:r}
It is natural to start with rhythm, since in the work of 
Messiaen,  this notion occupies a central place.  
  In \cite{Samuel}, p. 101, 
Messiaen  says: ``I consider rhythm as a fundamental element,  and may be the
essential element of music. It has conceivably    existed
before melody and harmony, and finally, I have a 
 preference for
this element." His monumental theoretical work, on which he worked for more than 
40
years, is titled ``Trait\'e de rythme, etc.". Since  the
first pages of this treatise, Messiaen rises up against the common opinion
which says that music is made out of {\it sounds}. He  writes: ``I say no ! No, not only with sounds ... Melody cannot exist  without
Rhythm!... Music is made first of all  with {\it durations}, {\it impulses},
 {\it rests},
{\it accents}, {\it intensities}, {\it attacks} and {\it timbres}, all things
which can be regrouped under a general term: {\it rhythm}". (\cite{Messiaen-Traite}, Tome 1 p. 40.)

The notion of rhythm is not foreign to mathematics, at least 
because elements like durations,  intensities and densities are measured with
numbers.  Timbre, decomposed into fundamental
frequency and harmonics, can also be expressed by a sequence of numbers. But of course, the use of number is not the main point, and there
are more profound reasons for the relation between rhythm and mathematics, and we shall discuss them below.

We shall mainly talk about rhythm
as a sequence of durations. Messiaen describes such a sequence
 as a
  ``chopping-up of	 
Time".\footnote{The relation between rhythm and time is important in 
Messiaen's thinking. Volume 1 of his Treatise \cite{Messiaen-Traite}, which is devoted to Rhythm, starts with a long chapter on
Time in all its aspects: absolute time, relative time, biological time, cosmological time,
physiological time, psychological time, the relation of time to eternity, etc., with long digressions on
the concept of time in mythologies, in the Bible, in Catholic theology,  in the theories of
 Einstein,  of Bergson
and in many other settings. See \cite{Shenton} for an interesting study on Time in the work of Messiaen.}
 Some of his compositions are based on 
 particularly simple (but never monotonic) rhythmic 
formulae.
We can mention in this respect Piece No. VII of his {\it Livre d'Orgue},
 titled {\it
Soixante quatre dur\'ees} (sixty-four durations), a composition which from an abstract point of view
resembles a game in which the composer takes different 
 groupings of elements 
in a chromatic\footnote{The word \emph{chromatic}, here and below, means that the values in the sequence increase linearly, that is, the sequence behaves like an arithmetic sequence. The terminology chromatic is used in music theory.} sequence of 64 durations   
  and  intertwines them using a geometric process which Messiaen calls ``symmetrical permutations", 
  which
we shall describe below. Let us mention also Piece No. V of his {\it Messe de
la Pentec\^ote}, titled {\it Le vent de l'Esprit}, in which  two chromatic
sequences are superimposed, the first one decreasing  gradually 
 from a  value equal to  23 sixteenths notes until the value of one sixteenth
note, and the second increasing gradually from 4 sixteenth notes
until 25 sixteenth notes.

Creativity often implies a profound immersion in the sources, 
and  
Messiaen's sources, for what concerns his rhythmical   language, are India and Ancient 
Greece.\footnote{One
has also to mention birds, but this is another story.}
 Indian and Greek rhythms have a dominant position, both
 in his
compositions and in his  teaching. Furthermore, in his written work, Messiaen keeps fixing  one's attention  on the 
arithmetical properties of these rhythms. We shall give some examples.

\subsection{Greek rhythms}

It is well-known that in Ancient Greece,  music was used as an 
accompaniment to  
poetry and to theatre, and therefore, musical rhythm  followed the rhythm 
of declamation. In this setting, there are essentially
two sorts of durations for musical notes, {\it long} durations, all
equal in value, corresponding to long syllables, and {\it short} durations, also all
equal in value,
 corresponding to short syllables, the value of a short duration being half
of the value of a long one. A rhythm in this sense, that is, a string  of long and short durations, is  called a  
 {\it meter}. Mathematicians know that there is a rich theory of combinatorics of strings of words written in an alphabet of two letters. 
 
 Ancient Greek music contained a rich variety of  meters, and 
  these   were classified in particular by
Aristoxenus in his impressive \emph{Harmonic Elements} which we already mentioned, see \cite{Barker} and \cite{Belis}, Vol. II. One of the characteristics of this music is that within the same piece,
meters    are 
of variable length,
in contrast with the meters of   (pre-twentieth century) Western classical music, where
a piece is divided into bars within which  the number of beats  is
constant.  Meters of variable length existed even in
Gregorian chant, which in some sense is a heir of ancient Greek music, and at some point during the Renaissance period,   Greek meters 
were in fashion.\footnote{For instance, at
 the beginning of the seventeenth century, Claude le Jeune composed choral works whose
 rhythm followed the principle of Greek  meters,   which is not based
    on the sole count of
 syllables, but  which takes into account their length or shortness.}
But  then the interest in them disappeared 
again, although there are
 reminiscences  of   Greek meters 
 in Romanian folk music and in compositions by Ravel and Stravinsky.
For instance, in Stravinsky's
{\it Rite of the Spring},  at the beginning
of the {\it Introduction}, the meter  switches 
constantly between the values 4:4, 3:4 and 2:4. Likewise, in the last piece,
{\it Sacrificial Dance}, the meter    changes constantly, taking values  like 5:16,
3:16, 4:16, 2:8, 3:8, 3:4, 5:4, and there are others. Messiaen revived the systematic usage of meters of variable lengths, 
teaching their principle
in his class at the Conservatory of Paris, 
and putting them into practice in his
compositions.  The first volume of his
Trait\'e \cite{Messiaen-Traite} contains a 170 pages chapter on Greek  meters.

 These ``a-metrical rhythms"  were used by
  Messiaen since his earliest
compositions. It seems that he cherished 
this  kind of freedom
in rhythm, and one reason for that is that it excludes monotony. 
 Messiaen, who sometimes described  himself as a {\it Rhythmician}, 
says in  \cite{Samuel},
p. 102, that ``a rhythmical music is a  
music which excludes repetition 
  and equal
divisions and which finds its inspiration  in the movements of nature, 
which are 
movements with free and non-equal durations." 
On p. 103 of the same treatise, he gives examples 
of a non-rhythmical music: ``Military music is the negation of
rhythm", and he notes  that military marches are most unnatural.
Likewise, there is no rhythm, he says, in a Concerto by Prokofiev, because
of the monotonicity of the  meter. On the other hand, he considers Mozart 
and Debussy as true rhythmicians. To understand this, we refer the reader to the
chapter titled ``\`A la recherche du rythme" in \cite{Samuel}. The reader might remember that the word rhythm refers here to a variety of notions: sequences of durations, but also of attacks, intensities,  
timbre, etc. 
In the first volume of his Trait\'e, Messiaen writes that rhythm contains periodicity,
``but the true periodicity, the one of the waves of the sea, which is
 the
opposite of pure and simple repetition. Each wave is different from the
preceding one and from the following one by its volume, its height, its duration, its
slowness, the briefness of its formation, the 
power of its climax, the
prolongation of its fall, of its flow, of its scattering..." (\cite{Messiaen-Traite}, Tome 1, p. 42).

Another aspect of Greek
 meters, which was seldom used  in Western classical music before
Messiaen, is the systematic use of rhythmical patterns whose value is a
prime number (other than 3), for instance 5, 7, or 17. 
One example of a
rhythm whose total value is 5 is the {\it Cretic rhythm}, defined by the
sequence 2, 1, 2 (that is, a rhythm corresponding to a long, then a short, and
then a  long syllable), and its two 
permutations,  2, 2, 1 and 1, 2, 2. The rhythm 2, 1, 2 is called {\it
amphimacer}, meaning (as Messiaen explains)   ``longs surrounding  the
short".  This  introduces us directly to two important notions in the
rhythmical language  of Messiaen.  The first one is related to the central 
symmetry of the sequence 2,1,2, which makes it an instance of  
  a {\it non-retrogradable rhythm}, and the second one is that of a
{\it permutation} applied to a rhythm. But before dwelling on that, let us say a few words on Indian rhythms,   which     also possess some beautiful properties.

 \subsection{Indian rhythms}
 
A significant characteristic of Indian music is the important place that it makes
for   percussion instruments
 like drums, cymbals,  bells, hand-clapping, and so on, and this makes
  rhythm  a very important factor in that music.
Let us quote Messiaen again: 
``Indian music is the music which   certainly 	
went farther than any other music
  in the
domain of rhythm, especially in the quantitative domain (combinations of
long and of short durations).  The  Indian rhythms, of unequalled 
 refinement and subtlety, leave far behind them our poor western rhythms
with their isochroneous bars, and their perpetual
divisions and multiplications by 2 (sometimes by 3)." (\cite{Messiaen-Traite} Tome 1, p.
258).

In the same way as do Greek rhythms, Indian rhythms abound in Messiaen's
compositions.  In his Trait\'e \cite{Messiaen-Traite},  
 the chapter concerning Indian rhythms occupies  
130 pages. In this chapter, Messiaen draws up  
lists of the 120 {\it de\c ci-t\^alas},\footnote{These rhythms have been
classified by the 13th-century Indian musicologist C\^arngadeva.
 Messiaen explains that, 
 in Hindi, de\c ci means
rhythm, and t\^ala means province. Thus, the word de\c ci-t\^ala refer the rhythms of
the various provinces. There are other interpretations for t\^ala; see for instance the
article {\it India} in the New Grove Dictionary of Music and Musicians.} of the
36 rhythms of the Carnatic (that is, South-Indian) tradition and 
 of other groups of Indian
rhythms, and he comments them thoroughly. 
 Here also,    a rhythm  is a sequence of
 numbers, and
 Messiaen expresses a real fascination for
the arithmetical properties of these sequences, a fascination
which he    transmits to his reader. Let us see a few
examples of these properties. He points out, whenever this is the case, that the
 sum of all the durations of some rhythms is
 a prime number. For instance, he records several de\c ci-t\^alas  
  whose total value are 5, 7, 11, 17, 19, 37, and so on. 
  This insistence on prime numbers may be surprising, but 
  we have already
mentioned the importance of these numbers for Messiaen. In the first volume of his Trait\'e, he writes that  ``the
impossibility of dividing a prime number (other than by itself and by
one) grants it a sort of force which is very effective in
the domain of rhythm."  (\cite{Messiaen-Traite} Tome 1, p.
266).

Another special class of   de\c ci-t\^alas 
which is highlighted by Messiaen  is the class of rhythms consisting of 
a sequence of durations which is
followed by its {\it augmentation}. For instance, the rhythm 1, 1, 1, 2,
2, 2 (de\c ci-t\^ala No. 73) is made out of  the sequence 1, 1, 1 followed by its
augmentation by multiplication by 2. 
An analogous feature occurs in de\c ci-t\^ala
No. 115, which is the rhythm  4, 4, 2, 2, 1, 1, constituted
by the sequence 4, 4, followed by its diminution 2, 2, and then by the
diminution of its diminution, 1, 1. Augmentation and diminution are
arithmetical transformations which are important in the art of
counterpoint, 
 which is the
art of transforming and combining musical lines,
and    which we shall discuss below in more detail.

There are more complex
combinations. For instance, in the rhythm 1, 3, 2, 3, 3, 3, 2, 3, 1, 3 (de\c ci-t\^ala
No. 27),  Messiaen notes that the odd-order durations are all equal, whereas  the even-order
durations consist in a regularly increasing and then regularly decreasing
sequence. He points out  that this rhythm was
used by Stravinsky in the {\it Rite of the Spring}, and that it is at the
basis of his theory of {\it Rhythmic characters}.\footnote{There are three rhythmic
characters here:  one character stays still, 
another one 
is decreasing and the third one is increasing.} Messiaen makes extensive 
use of 
  rhythmic
characters in his {\it Turangal\^\i la Symphony} (composed in 1946-1948). 

Finally, let us mention that the de\c ci-t\^alas contain  
several instances of  non-retrogradable rhythms,
that is, rhythms    consisting of a sequence of durations followed by its
mirror image (with sometimes a common central value).  We already
encountered such rhythms  when we talked
about Greek rhythms. For instance, de\c ci-t\^ala No. 58 is
 the Greek amphimacer rhythm,
2, 1, 2  , which Messiaen describes as ``the simplest and
the most natural non-retrogradable rhythm, because it is based on the
number 5, the number of fingers in the hand." (\cite{Messiaen-Traite}, Tome 1, p. 289). 
Other examples
of
non-retrogradable  de\c ci-t\^alas are  2, 2,
1, 1, 2, 2 (de\c ci-t\^ala No. 26), 1, 1, 2, 2, 1, 1 (de\c ci-t\^ala No. 80) and 2, 1, 1,
1, 2 (de\c ci-t\^ala No. 111), and there are several others. In the next section, we
 discuss  
non-retrogradable rhythms at fuller length.

\subsection{Non-retrogradable rhythms}

The use of non-retrogradable rhythms is an example of the systematic use of symmetry in the music of Messiaen.

 Messiaen dealt with non-retrogradable rhythms since his early
compositions, and he attached great importance to them in his
first theoretical essay, {\it Technique de mon Langage Musical} (1944). A
non-retrogradable rhythm is a sequence of durations which gives the same result
whether it is read from left to right or from right to left. It may be
good to recall here that {\it retrogradation} is a classical device in the art
   counterpoint which we shall consider more thoroughly later on in this article. It    transforms a  certain musical motive by
reading it backwards, that is, beginning from the last note and ending
with the first note. The initial motive is then called the {\it motive in
direct motion}, and the transformed motive the {\it motive in retrograde
motion}. Thus, a non-retrogradable rhythm can be regarded as the
juxtaposition of a 
 motive in direct motion  and of a motive in retrograde motion,
with sometimes  a central value in common. 

Retrogradation, as a counterpoint operation, was used and taught since
the beginning of this art, around the fourteenth century. But before Messiaen, it was
usually applied to a
melodic motive, that is, to a sequence of pitches, whereas with
 Messiaen, retrogradation acquired a more abstract
character,  since he applied it systematically to rhythm, regardless of
pitch. Thus, the listener of Messiaen's music is invited to feel  
retrogradation at the level of durations only, since there  need not be any
regular correspondence (transposition, symmetry, etc.) between the pitches  of
the motive in direct motion and those in the motive in retrograde
motion.

For instance, in the ``Danse de la fureur, pour les sept trompettes" (Part VI
of Messiaen's {\it Quatuor pour la fin du Temps}),  
we find the following succession of non-retrogradable rhythms:

$$3, \ 5, \  8, \ 5,\ 3$$
$$4, \ 3, \  7, \ 3,\ 4$$
$$2, \ 2, \ 3, \ 5, \  3, \ 2,\ 2$$
$$1, \ 1, \  3, \ 2,\ 2,  1, \ 2, \ 2, \ 3, \ 1, \ 1$$
$$2, \ 1, \  1, \ 1,\ 3, \ 1, \ 1, \ 1, \ 2$$
$$2, \ 1, \  1, \ 1,\ 3, \ 1, \ 1, \ 1, \ 2$$
$$1, \ 1, \ 1, \ 1, \ 1, \ 3, \ 1, \  1, \ 1,\ 1, \ 1$$
$$3, \ 5, \  8, \ 5,\ 3$$
(the   unit being the sixteenth note). In the analysis of this sequence
of bars that he makes in the second volume of his Trait\'e, Messiaen, of course,
highlights the fact that the
total number of durations in bars 3 and 4 is 19, and that in bars
5, 6 and 7, this value is 13, pointing out  again that 19 and 13 are
prime numbers (\cite{Messiaen-Traite},  Tome 2, p. 26).
It would be  tedious and superfluous to try to draw up a long list of such examples,
since  the are plenty of them in Messiaen's compositions. But it is 
natural 
to
raise the question 
of  why non-retrogradable rhythms are interesting. One may as well ask 
 why is the symmetry of a face a beauty criterion. 
 Messiaen answers this question on rhythm in his own way, and he gives 
  two kinds of reasons, one
of an aesthetic nature, and the second one philosophical. In \cite{Messiaen1944}, he talks about
the {\it charm} which a non-retrogradable rhythm  produces on the listener of
his music. He considers that this charm is of the same order
as the one which is produced by his {\it modes of limited
transposition} (which we shall discuss below), and   he calls
the non-retrogradable rhythms and the modes of
limited transposition as two {\it mathematical impossibilities}.
The impossibility, in the first case,  lies in
 the fact that it is  ``impossible to reverse such a rhythm", since when
we reverse it, we obtain exactly the same rhythm.
Let us quote Messiaen from the first first of his  \emph{Technique de mon langage musical} (\cite{Messiaen1944}, Tome 1, p. 5): 
\begin{quote}\small
One point will attract first our attention: the
{\it charm of impossibilities}... This charm, at the same time voluptuous
and contemplative, lies  particularly in certain mathematical
impossibilities in the modal and of the rhythmic domains. The modes which
cannot be transposed beyond a certain number of transpositions, because 
if one does
so, he
falls again on the same notes; the rhythms which cannot be
retrograded because if one does so, he recovers 
  the same order of the values...
  \end{quote}
On page 13 of the same Trait\'e, Messiaen describes the impressions which
these impossibilities    produce on their listener. 
\begin{quote}\small
Let us consider
 now the  listener of our modal and rhythmic music; there is no
 time for him, at the concert, to check  
non-transposition  and   non-retrogradation, and at that moment, these
questions will not interest him any more: to be seduced, this will be
his unique desire. And this is precisely what will happen: he will
undergo despite his will the strange charm of these impossibilities; a
certain ubiquitous tonal effect of the non-transposition, a certain
unity of movement (where beginning and end merge,  because they are
identical) of non-retrogradation, all things which will certainly lead
him to that sort of ``theological rainbow" which our musical language
tries to be, a language which we are trying to edify and to
theorize.\footnote{This sentence is reminiscent of the sentence by Leibniz that we quoted in \S 2, saying roughly that music is a secret exercise in arithmetic.}
\end{quote}

Messiaen after this explains the philosophic relevance of
 non-retrogradable rhythms,  which justifies also the term ``theological
rainbow" in the last sentence. 
A non-retrogradable rhythm, according to him, can give the listener
 a feeling  of infinity. Indeed,
whether one reads it from left to
right or from right to left, a non-retrogradable rhythm stays invariably
the same, and  in this sense,  such a rhythm  has no beginning and no
end.   
In the second volume of his Trait\'e,  Messiaen says that a 
non-retrogradable rhythm 
draws its strength from the fact that ``like
 Time, a non-retrogradable rhythm is irreversible. It cannot move
backwards, unless it repeats itself... The future and the past
are mirror images of each other" (\cite{Messiaen-Traite}, Tome 2, p. 8).

In conclusion to this section on rhythm, we quote again Messiaen, who mentions in the first volume of his Trait\'e rhythms which are ``thought for the only intellectual
pleasure of number" (\cite{Messiaen-Traite}
Tome 1, p. 51). The reference to number reminds us again of ideas that Leibniz, Euler, Diderot, and others emitted about music and which we recalled in \S \ref{s:brief}. Numbers  can {\it a priori}
 appear as being  severe, austere, 
and devoid of lyricism. Expressed as rhythm, they are given a new
 life.

Messiaen notes in his Trait\'e \cite{Messiaen-Traite} that nature is full of non-retrogradable rhythms, starting with the human face, with the two ears, the two eyes and the nose at the center, or like the teeth inside the mouth. He also makes an analogy with several architectural edifices, including the marble bridge in Beijing that leads to the Summer Palace of the Chinese emperors, together with its reflection in the water.

Messiaen certainly compared the beauty
of certain rhythms, which are built   as   sequences of
numbers with  rigorous properties, 
to the beauty of certain faces with regular and
symmetrical features, 
  to that of French gardens that follow completely
symmetrical plans, to that of Romanesque cathedrals, and to that of the wings of
butterflies. 

As a last wink  to non-retrogradable rhythms, let us mention that the
number of pieces in the seven books of Messiaen's composition {\it Catalogue d'Oiseaux} is
respectively 3, 1, 2, 1, 2, 1, 3.
\bigskip

\subsection{Symmetrical permutations}

 Permutations of finite sets play an important role in the
music of Messiaen, and groups are present there. We shall see this more precisely below. 

Given a sequence of musical   objects (e.g. pitches, dynamics,  durations, etc.),
one   obtains another
sequence by applying to it a permutation.  For instance,  
  retrogradation is a special kind of permutation.
The  problem is that as soon as the number of objects is large, the total number of permutations becomes too large to play a significant role in music. For instance, for a
sequence of 5 objects, there are 120 distinct permutations, for 6 objects, there are  720 distinct permutations, for 7 objects, there are 5040 distinct permutations, and then the number of permutations become huge. Thus, for the use of permutations in music, one has to make
choices, because if the order of the symmetries used is too large, the ear cannot discern them. This leads us to the theory that Messiaen calls \emph{symmetrical} permutations, that is, permutations which have a small group of symmetries.
In his compositions, he applies to a musical motive the
iterates of a given  symmetrical permutation. For instance,  
retrogradation is    
 of order two and therefore it is symmetrical.  Symmetrical
permutations which are more complicated  than retrogradation are used 
 already in 
his early pieces, for instance in  the
{\it Vingt Regards sur l'Enfant J\'esus},    composed
in 1944.

The piece {\it Chronochromie} (1960) starts   with a 
chromatically increasing sequence  of 32 durations, starting with a 
thirty-secondth note, and ending with a 32 $\times$ thirty-secondth note, that is, a
whole note.  
Messiaen applies to it the  permutation 3, 28, 5, 30, 7, 32, 26, 2, 25, 1, 8, 24, 9, 23, 16, 17, 18, 22, 21, 19,
20, 4, 31, 6, 29, 10, 27, 11, 15, 14, 12, 13. He then applies the same permutation to the new sequence, and so forth. After 35 steps, we recover the initial
sequence, 1, 2, 3, ..., 32. 

In {\it Chronochromie}, Messiaen uses,  as rhythmical motives, a collection of
rhythms which belong to this set of iterates.   

It is natural to ask why these permutations are important in music.
 In his treatise, Messiaen describes the symmetrical permutations as
a third {\it mathematical impossibility}. The listener of such permutations in a musical piece is supposed to be 
dazzled by the same {\it
charm} as with non-retrogradable rhythms and of modes of limited
transposition, the two impossibilities that we already mentioned. In \cite{Samuel}, p. 222, he describes the piece {\it
Chronochromie} as
``durations and permutations rendered sensible by sonorous colors; this
is indeed a Color of time, a Chonochromy".

For mathematicians, it is amusing to see a large number of explicit examples contained in that musical treatise,
and to imagine how tedious it was for Messiaen to write them.

There is a special kind of permutations to which Messiaen attaches more
importance (and which probably are at the origin of the adjective ``symmetrical"). To obtain such a permutation, one starts at the center of a
sequence of objects, and then takes successively  one
object from the right and one object from the left, 
until one reaches the
two ends of the sequence, and one of its iterates. 
For instance, this
process transforms the sequence 1, 2, 3, into the
sequence  2, 1, 3 . Applying the same rule
to 2, 1, 3, we obtain 1, 2, 3, which is the sequence we started with. Thus,
the   permutation $1,2,3\mapsto 2,1,3$ is of order 2. Let us
consider now a sequence of four objects. The iterates are: 
$$1, 2, 3, 4 \mapsto 2, 3, 1, 4\mapsto 3, 1, 2, 4\mapsto 1, 2, 3, 4.$$
and in this
case, the order of the permutation   is 3. This is a way of finding ``by hand" symmetrical permutation. It is interesting to see in which terms Messiaen describes this 
 mathematical process.
In \cite{Samuel}, p. 119, he 
 says: 
 \begin{quote}\small There are durations which follow one another
 in a
certain order, and we read them again in the initial order. Let us
take, for example, a chromatic scale of 32 durations: we invert them
according to a given order, then we read the result according to this
 order, and so forth  until we recover  the initial chromatic
scale of 32 durations. This system produces interesting and very strange
rhythms, but above all, it has the advantage of avoiding an absolutely
fabulous number of permutations. You know that with the number 12,
so much beloved by serialists, the number of permutations is
479.001.600 ! One needs years to write them all. Whereas with my
procedure, one can, starting  with larger numbers --  for instance 32 or
64 -- obtain better permutations, suppress the secondary permutations
which lead only to repetitions, and work with a reasonable number of
permutations, not too far from the number we started with.
\end{quote}

\section{Counterpoint}\label{s:c}

 Counterpoint is the classical art of transforming,  combining, and
superimposing musical lines.\footnote{The word counterpoint comes from the latin expression \emph{punctus contra punctum} which means ``point against point", expressing the fact that on a musical score, the dots that represent different pitches that are played at the same time, as the result of the superimposition of the musical lines, appear vertically one above the other.} Symmetry is extensively applied in this art, and we shall see this in this section. Two important treatises on counterpoint are those of Tinctoris (1477)  \cite{Tinctoris} and Fux (1725) \cite{Fux}. 

\subsection{The use of integers mod 12}\label{use}
  There are a few standard operations in counterpoint, and
it is practical to describe them using the language of integers mod 12. We
briefly discuss this formalism here, and this will also serve us in the description of Messiaen's {\it modes of limited transposition}. We note however that  
 Messiaen did not use the notation of integers modulo 12; in fact, he did not have any background in mathematics. He had his own words to define mathematical objects and to explain what he wanted to do with them. This usually involved a heavy language.

In two famous articles, published at the beginning of the 1960s (see \cite{Babbitt1960} and
\cite{Babbitt1961}), 
Babbitt applied the language and 
  techniques of group theory in music (in particular in
twelve-tone music). He used in particular the concept of
{\it pitch-class} which became an important tool in the teaching of certain
musical theories; for instance, this is an important factor in the  textbook by
Allen Forte \cite{Forte}, which was for many years one of the main references in the USA for twentieth-century musical analysis  techniques.

The principle of the use of integers mod 12 is the identification of notes whose pitches differ by an octave. This is a natural identification, because in practice, it is observed that in general,  the notes produced by the voices of a man and of woman (or of a man and a child) singing the same song, differ by an octave, although they are considered to be the same notes. This octave identification is also suggested by the fact that the names of notes  that differ by an octave have the same name, and therefore there is a cyclic repetition in the names of notes that are played by a traditional instrument. In fact, it can be difficult for a non-expert to say whether two notes with the same name played by different instruments (say a flute and a violin) correspond to the same pitch or to pitches that differ by an octave. 
This fact has been pointed out and analyzed by several music theorists, in particular in the set of problems concerning music theory which are attributed to Aristotle (\cite{Aristotle-Problems}, Volume I).\footnote{Music was not the main subject of study of Aristotle, and in fact, Aristotle had no preferred subject of study; he was industrious in all sciences -- and music theory was one of them --  working on them, teaching them, and writing essays on all them. He had a set of \emph{Problems} that he made available to the students of his school (which was called the \emph{Lyceum}), which one could compare to the lists of open problems that are known to mathematicians today, except that Aristotle's problems  concern all sciences. The set of problems from the school of Aristotle that reached us contains about 900 problems, classified into 38 sections, of which two are dedicated to music (one section is more directed towards acoustics, and the other one towards music theory). Aristotle, like his teacher Plato, payed careful attention to music and acoustics, and Aristoxenus, who became later on the greatest Greek musicologist, was Aristotle's student.} It has also been formalized as a principle, the ``principe de l'identit\'e des octaves", by Rameau in his \emph{Trait\'e de l'harmonie r\'eduite \`a ses principes naturels} \cite{Rameau-Traite}, and this principle has been thoroughly used before and after Rameau.
We note by the way that there was a correspondence between Euler and Rameau about this principle, and some of these letters reached us, see \cite{Jacobi}, \cite{FN} and \cite{Euler-Hermann}.

In any case, we are considering now the so-called \emph{equally tempered chromatic
scale}.\footnote{The equally tempered scale became more or less universally adopted in the nineteenth century for various reasons, one of them being the appearance of large orchestral ensembles, in which several  kinds of different instruments had to be tuned in a uniform way, and the most convenient way appeared to be the one based on equal division. Furthermore, with that scale, transpositions become easy to perform, and playing a given piece starting at any note was possible. We note by the way that equally tempered scales were already considered in the theoretical work of Aristoxenus.}  In other words, we are considering notes which correspond to a division of an octave into twelve equal intervals. Since the ratio of the frequencies of two pitches separated by an octave is equal to 2/1, the ratio of two successive notes in the equally tempered scale is equal to $\sqrt[12]{2}$. The notes in this scale correspond to the sequence {\it C, C$\sharp$, D, D$\sharp$, E, F,
F$\sharp$, G, G$\sharp$, A, A$\sharp$, B} that appear in that order within an octave on a piano keyboard. Represented  by
the integers 0, 1, ..., 11 in that order, these numbers are considered as  elements of the
cyclic group $\mathbb{Z}_{12}=\mathbb{Z}/12\mathbb{Z}$.
For our needs here, the applicability of symmetry and group theory is easier if we use the equally tempered scale, and these mathematical notions naturally appear in the description of the modes of limited transposition as we shall see below.

In the language introduced by Babbitt, and developed by  Forte and others, these integers  
 are called {\it pitch-classes}, and they
represent equivalence classes of pitches, that is, pitches defined up to addition 
of a multiple of an
octave.
 A {\it pitch-class set} denotes then a set of
pitch-classes, and it is represented  by a sequence of distinct elements
of the group $\mathbb{Z}_{12}$. To denote a pitch-class set, Forte uses  square brackets,
 for instance [0, 1, 3]. One writes, to be brief, {\it pc} and {\it pc-set} for pitch-class and pitch-class set.

 It is fair to note here that this notion of pc set  was known (without the name) in the nineteenth century, for instance by Heinrich Vincent  (1819-1901) \cite{Vincent} and Anatole Loquin (1834-1903) \cite{Loquin}. In the twentieth century, Edmond Cost\`ere \cite{Costere} made an exhaustive list of pc sets long before Babbitt and Forte, using the name \emph{\'echelonnement}, but without taking into account the notion of inversion.\footnote{I learned this from F. Jedrzejewski.}
  
  Let us now review the relation with counterpoint.
  
The three basic operations of counterpoint 
 are \emph{transposition, inversion} and
\emph{retrogradation}, and they have been identified  and highlighted 
early at the beginning  of this art (that is, at the beginning of  the Renaissance). 
Using a mathematical language, the operations are represented respectively by a translation, a symmetry with respect to a horizontal line and a symmetry with respect to a vertical line (which is outside the motive): equivalently, for the last operation, one reads the motive backwards.  

At the beginning of the twentieth century,
 and notably under the influence of Arnold Sch\"onberg
  and his School (the so-called
 Second Viennese school),
the counterpoint operations were used in a more abstract
and  systematic way by composers who became known as  \emph{dodecaphonic} or
 \emph{twelve-tone} composers and  later on as \emph{serialists}.
 In this setting, a 
{\it series}, also called a {\it tone row}, is 
 a pitch-class set, that is, 
 a collection of distinct notes,  
 with no special melodic value\footnote{This  contrasts with the usual motives
  of old contrapuntal writing where the building block, called the \emph{theme}, has an intrinsic musical value; it was sometimes -- e.g. in the music of Bach -- a theme extracted from a popular song or of a well-known melody. It should be noted however that Bach also wrote magnificent chorals based on poor themes, in fact, themes obtained by concatenating notes in a manner which \emph{a priori} is uninteresting. The richness of the resulting harmonies is due to the cleverness of the assembly of the theme with its symmetric images, obtained through the counterpoint operations. This is another way of proving -- if a proof is needed -- that it is the mathematical structures that are behind a musical piece that makes its beauty.}
and which, together with its
transformations by transposition, retrogradation and inversion, 
is used as a building block for a musical composition.
A {\it twelve-tone series} is a tone row in which every  
 integer  modulo 12
appears. Thus, a twelve-tone series is simply a  permutation 
of the sequence of integers 0, 1, 2, ..., 11.
The three operations of  transposition, inversion and
retrogradation can be expressed
  simply and
elegantly using the formalism of pitch classes, and this was done in the work of Babbitt and Forte. Indeed, transposition
corresponds to  translation modulo 12 in the group $\mathbb{Z}_{12}$, inversion
is the map $x\mapsto -x$ mod 12, followed by a translation, and 
retrogradation consists in replacing a certain motive $n_0, n_1,...,n_k$
by the same motive written backwards, $n_k,...,n_1,n_0$.  We note that in
the context of serial music, the integers $n_i$ in the sequence representing a melodic motive have to be distinct as classes mod 12. It is appropriate to quote here Forte (who is not a mathematician),
 who on p. 8 of his essay \cite{Forte},
stresses the fact that this formalism is not
useless: ``The notion of mapping is more than a convenience in
describing relations between pitch-class sets. It permits the development
of economical and precise descriptions which cannot be obtained using
conventional musical terms".
Let us note that this kind of theory is known among musicologists as \emph{set theory} (and in France, it is called \emph{American set theory}, or simply ``set theory", using the English words), although it has not much to do with set theory as mathematicians intend it. A recent reference is \cite{AB}.

Let us now return to Messiaen.

\subsection{Generalized series}

 Messiaen never belonged to the serial school
and in fact, he belonged to
no school. However, it is true that some of his compositions contain
techniques which appertain to that school, and above
all, his piece  {\it Mode 
de valeurs et d'intensit\'es} (Piece No. II of his   {\it Quatre Etudes de Rythme} for
the piano, 1949). In fact, in this piece, Messiaen  
 goes beyond the existing twelve-tone techniques by using the so-called
 {\it generalized series} (Messiaen used the French word \emph{supers\'erie}), that is, 
    not only series of
pitches, but also series of rhythms, of intensities, of attacks, 
 of dynamics, and of timbres.
An example of a series of intensities is the ordered sequence
  $$ppp, \ pp, \ p, \ mf, \ f,  \ ff, \ fff.$$

The  piece
{\it Mode  de valeurs et d'intensit\'es} had a significant  influence on 
the so-called post-serialist (or ``integral serialism") school,
to which belonged at some point
 Boulez 
and Stockhausen, who  had been students
of Messiaen. Boulez, stimulated by this piece, used extensively the principle of   generalized series  applied to timbre, intensity, duration,
and so on. For instance, he used in his piece 
{\it Polyphonie X pour 18 instruments} (1951) a series of 24 quarter-tones. In
his piece {\it Structures I} for two pianos (1952), he used    series of 12
durations, of 12 intensities ({\it pppp, ppp, pp, etc.}) and 
of 12 attacks ({\it staccato, legato, etc.}).
 The position of Messiaen 
regarding the serial movement was  moderate, and in fact, he
disliked the excessiveness of abstract formalism that this movement gave rise to. But most of all, 
he disliked the absence of colors in that music. He described  the music of
Sch\"onberg and the Second Viennese School as black, morbid and as a  ``night music".
We mention also that Boulez later on broke with the serial 
school, and in fact, one should also note that  
the period which followed 
the composition of 
{\it Mode 
de valeurs et d'intensit\'es} was also a period of profound
questioning  of the usefulness of the serial techniques, 
by the serial composers themselves.

The new combinatorics that Messiaen introduced in his
piece {\it Mode 
de valeurs et d'intensit\'es} had a great  
  influence that Messiaen himself  disapproved. 
He declares in
\cite{Samuel} p. 119: 
\begin{quote}\small
In  {\it Modes de valeurs et d'intensit\'es}, I  used a
superseries  in which sounds of the same name come past various
regions making them change as to octave, attack, intensity, duration. I
think that this was an interesting discovery... Everybody used to talk
only about the superserial aspect,
\end{quote}
 and also in \cite{Samuel}, p. 72: 
\begin{quote}\small I was very annoyed by the absolutely inordinate importance 
which has been granted  to that minor work, which is only three pages long
and which is called {\it Mode de valeurs et d'intensit\'es},  
under the pretext that it was supposed to be at the origin of the serial
shattering in the domain of attacks, of durations, of intensities, of
timbre,  in short, of all the musical parameters. This music may have been
prophetic, historically important, but musically, it is nothing...
\end{quote}

\subsection{Rhythmical counterpoint}

Messiaen developed in his teaching a theory of counterpoint which is
proper to rhythm.   One of the features of this theory is that the counterpoint transformations are applied to a musical 
motive at the level of
rhythm, regardless of pitch. Counterpoint writing
includes  other transformations than the
three that we mentioned in \S \ref{use}, and that they can all be 
applied to rhythm. In fact, we already encountered two of these transformations
  in the section on Indian rhythms,
{\it augmentation} and {\it diminution}. These
transformations consist in taking a certain melodic
motive, keeping unchanged the sequence of pitches,  and multiplying 
the values of
all  the
durations by a constant factor (which is $>1$ in the case of
augmentation and $<1$ in the case of diminution). These transformations exist in classical counterpoint, where the melodic motive stays the same while the durations are transformed. But in the music of
Messiaen, augmentation and diminution have a more abstract character, because they affect rhythm regardless of pitch; the latter can undergo either regular or irregular
transformations. For example, at the beginning of  Piece No. V ({\it Regard du Fils sur le Fils}) of the 
 {\it Vingt Regards sur l'Enfant J\'esus}, there is
   an augmentation of rhythm by a factor of
3/2, while  the motive and the transformed one are unrelated from the point of view of pitches. The themes are written in different
modes, whereas
their various combinations in the piece 
 constitute a {\it rhythmical canon}. In any case, augmentation and diminution preserve  
   non-retrogradability.

There are other rhythm transformations  which preserve 
non-retrogradablility, and which are discussed  by Messiaen in the second volume of his Trait\'e (\cite{Messiaen-Traite}, Tome
2, p. 41). These include {\it symmetrical amplification}
 and {\it symmetrical
elimination} of the extremities. A symmetrical amplification consists in
adding at the two extremities of a given rhythm  another 
rhythm together with its retrograded form in such a way that 
non-retrogradability
  is preserved. For instance, in Piece No.
XX of the  {\it Vingt Regards}, the theme   is
exposed at bar  No. 2,  it a very short theme, and
its rhythmic value is 2, 1, 2 (the amphimacer non-retrogradable rhythm, 
with the unit being
the sixteenth note). The theme is then symmetrically amplified  
 at bar No. 4, where it becomes 2, 2, 2, 1, 2, 2, 2. At bar No. 6,
it is amplified in a different manner, where it becomes 2, 3, 2, 2, 1, 2, 2, 3,
2. We find again this rhythm, with  two different 
amplifications, later on in the same piece (bar 82):
$$2, 1, 2$$
$$2, 2, 2, 1, 2, 2, 2$$
$$2, 3/2, 2, 2, 2, 1, 2, 2, 3/2, 2,$$
and so on.
There are other operations that preserve non-retrogradability, in particular the ones 
called by Messiaen {\it enlargement} and {\it diminution of the
central value} (\cite{Messiaen-Traite}, Tome 2, p. 30).

It is interesting to know that the {\it de\c ci-t\^alas} contain examples of all the rhythm transformations   that
we mentioned.

Classical counterpoint techniques are usually applied in the composition of
{\it canons}. Messiaen wrote {\it rhythmical canons},
that is, pieces based on superimposition of rhythms and their transformations,
following a certain regularity rule (the Greek word ``canon" means ``rule") -- for instance a periodicity
pattern -- governing the superimposition between  a motive and the transformed
motive.
A composition such as the {\it Vingt Regards sur l'Enfant
J\'esus} contains several canons of non-retrogradable rhythms, 
and it is interesting  that
Messiaen  indicates   explicitly  on  the score, at several places, the
subjects and the counter-subjects  of these canons, as well as
the various transformations which they undergo,
in order to 
 help  the reader understand
 the structure of the piece.

Let us note finally that Messiaen's {\it Turangal\^\i la Symphony} 
is considered in
itself as an immense counterpoint of rhythm.

In the next section, we shall see how counterpoint operations were
 necessary for Messiaen
in his transcription of   bird songs.

\subsection{Counterpoint and bird singing}

Birds represented one  the most important sources  of inspiration
   for  Messiaen, and probably the most important one.
   His interest for bird songs exceeded his interest for any other
   music. His early piece
   {\it Quatuor pour la fin du Temps} (1946) contains already 
   a lot of bird song
   material, and one can hear bird songs in all  the pieces that he
   composed after
1950.  
In fact, bird songs are even {\it the} central element
in several of his pieces, including {\it Le r\'eveil des oiseaux} (1953), 
{\it Oiseaux exotiques} (1956),
{\it Catalogue d'oiseaux} (1958), {\it La fauvette des jardins} (1970),
{\it Petites esquisses d'oiseaux} (1985), etc.
In this respect, one has to  mention also    his opera {\it Saint Fran\c cois d'Assise}
(composed in 1975-1983)
 and another major piece, {\it Des canyons aux
 \'etoiles} (composed in 1971-1974).

 Let us see how counterpoint operations were necessary during this transcription
  process.
  
 Diminution, first,
  was needed  because some birds sing at a tempo which is so fast
 that it is impossible to reproduce by any performer. Thus, Messiaen,
  in writing  these songs, 
 was led  to use rhythm
 diminution. Transposition was also needed because some birds sing at a register
 which is too high to be played on an instrument like the piano, for which, for instance,   
 {\it Catalogue d'oiseaux} was written.\footnote{The piano is among the instruments which have the widest register.}   
   Messiaen   was led therefore to transcribe the song at  a
 register which is one
    or several octaves
 lower. Finally, Messiaen appeals to an unusual operation of {\it augmentation}
 at the level of pitches. Indeed, 
 some  birds use in their singing   tiny
 intervals (of the order of a comma), and to transcribe their  songs  
  for an instrument like the piano,  Messiaen, in the absence 
   of these
 micro-tonal intervals,  replaced,
 for instance, an interval of  two commas by a semitone and the rest of the transcription follows, multiplying the intervals by the same factor. Thus, an interval of four
 commas 
   became  a
 second, and so forth.

\section{Modes of limited transposition}\label{s:ml}

The theory of modes of limited transposition is a spectacular example where symmetry and group theory are used in music. It is one of the most important characteristics of Messiaen's compositions (saying that a certain piece of Messiaen is in a certain mode is like saying  that Beethoven's fifth symphony is in the C minor scale), and it is a purely mathematical notion. Before talking about this subject, we say a few words about the meaning of the word {\it mode} in
music, which might be useful for the reader. In the last subsection, we shall talk about the relation between modes and colors as theorized by Messiaen. This subject has a long history, and it is in the tradition of the relation between pitch and color that was studied by Newton and others which we mentioned in \S \ref{s:brief} of this article. 

\subsection{Modes}
The title of Messiaen's treatise, \emph{Trait\'e de rythme, de couleur et
d'ornithologie} \cite{Messiaen-Traite}, contains the word color. He describes there his perception of the relation between mode and
color, and fully understanding this requires a certain mental
faculty which Messiaen possessed, namely, that of seeing in a precise way the different
colors
associated with  modes,  within harmonies or within melodic motives. We shall say more about this in the last subsection of this article, and the interested reader can consult \cite{Samuel}, p. 95, the
chapter entitled {\it Of sounds and colors}. In the absence of this faculty, one can
resort  to a more down-to-earth definition, and consider a mode as 
a sequence of distinct notes  which describe (or, in a certain sense,  characterize)  the {\it atmosphere} of a musical composition, or of the
part of the composition which is concerned by this mode. 
There also exist pieces that are {\it
polymodal}, and some others that are {\it not} modal. But most of the pieces that we hear are modal. The major and minor {\it
scales} which are usually  associated to the familiar pieces of music are approximations of this atmosphere. They classify the so-called
``tonal music", that is, most pre-twentieth century classical pieces (and practically all modern popular music pieces). Modes are described by sequences of notes. Rameau, in his {\it
D\'emonstration du principe de l'harmonie} (see \cite{Rameau-Traite}, Vol. 3),
 writes: ``In music, a mode
is nothing else than a prescribed order between sounds, both in harmony and in
melody." The music of 
Ancient Greece, that of China  and that
of India possess a rich variety of modes.
 Byzantine chant, and to some extent Gregorian chant,  that are
considered to have their roots in Ancient Greek music, have also several modes.

\subsection{Modes of limited transposition}

Mathematically, a
 {\it mode of limited transposition} is a sequence of
distinct notes that has a nontrivial group of symmetries by the action of translations mod 12, that is, by transpositions mod 12. 
Although Messiaen did not do so, 
it is convenient  to describe these modes using  the formalism of {\it
pitch-classes}  which we recalled   
 above.  In this setting, a mode is
  a certain {\it pc-set}, that is, a set of elements in the group
$\mathbb{Z}_{12}$, and a mode of limited transposition is associated to a
{\it pc-set} which is invariant by a nontrivial translation mod 12. In this sense, a mode of limited transposition is a pc-set which has some symmetry.
The order of a translation in $\mathbb{Z}_{12}$ is necessarily a divisor of 12, and we shall see now how the  various modes are associated to
translations of  order  2, 3 or 6. The familiar major and minor scales are invariant by a translation of order 12 and not less (and therefore by the identity map) and in that sense they are not modes of limited transposition.
 The expression of a given mode in
terms of pitch-classes, that we give below, 
 highlights in a particularly simple manner the symmetries of such a mode.
 
 Messiaen used modes of limited transposition
 since his first compositions. They are present, for instance, in his
   {\it Preludes for the Piano}, which he composed in 1929, at the age of 19. He writes in his \emph{Technique de mon langage musical}  \cite{Messiaen1944} p. 13: ``These modes realize in the vertical sense [that is, at the level of pitch] what non-retrogradable rhythms realize in the horizontal sense [at the level of time]". There is no non-retrogradability in the modes of limited transposition, and the common feature between the two notions is the presence of symmetry.

In his classification, which we shall recall below, Messiaen keeps  seven modes of limited transposition. In his compositions, he uses
mostly   Modes 2, 3, 4 and 6, with a preference for Mode 2, which at
his first
transposition corresponds to different gradations of blue-violet,
which was his favorite color.
 We also note that this classification in seven modes has not much to do with the previous existing classifications
 (e.g. the eight modes of Gregorian chant or of Byzantine music, or those of Ancient Greece which are much larger in number).
 
As we said above, the theory assumes that the tuning of the instruments for which the piece is written is the modern equally tempered tuning (equal intervals between two consecutive notes).

\medskip

\noindent {\it Mode 1.---} This is the sequence of notes 
{\it C, D, E, F$\sharp$, G$\sharp$, B$\flat$}.
Messiaen
describes  this mode as being ``two times transposable", and by this he
means that the corresponding {\it pc-set} is invariant by the translation
$x\mapsto x+2$. The {\it pc-set}  [0, 2, 4, 6, 8, 10] and it is the orbit of
0 by this translation.\footnote{The pitch-class sets that are
described
here start  at 0, and they correspond to
modes which are   {\it at their first
transpositions}.} In fact, as Messiaen notes, this mode is Debussy's \emph{whole tone scale} and it was already thoroughly used by this composer, for instance in certain passages of his opera \emph{Pell\'eas et M\'elisande}. It was also used by Paul Dukas (who was Messiaen's teacher) in his opera \emph{Ariane et Barbe-Bleue} and by Liszt in his opera-fantasy for piano \emph{R\'eminiscences de Don Juan}. This scale is called the \emph{whole tone scale} because each degree increases by a tone. It corresponds to the division of an octave into six equal subintervals, and it is represented in musical notation as follows:
\[
\includegraphics[width=.8\linewidth]{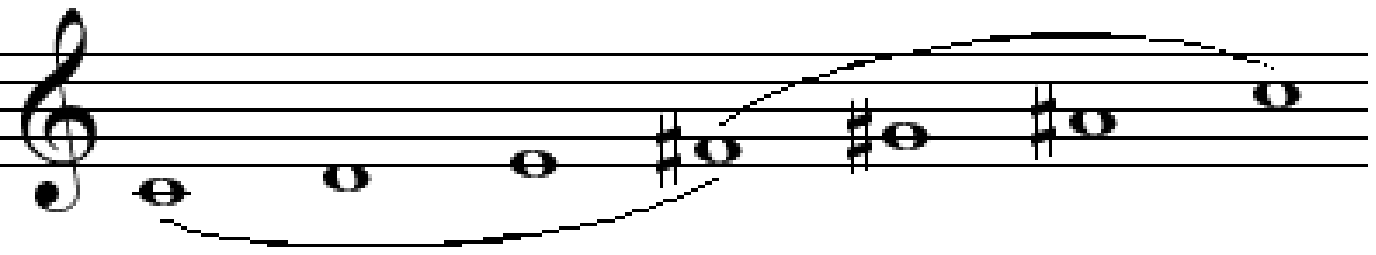}
\]
Messiaen made a very moderate use of the first mode, probably because it was not very innovative.\footnote{Messiaen writes, in cite{Messiaen1944} (Vol. I p. 52): ``Claude Debussy (in \emph{Pell\'eas et  M\'elisande}) and, after him, Paul Dukas (in \emph{Ariane et Barbe-Bleue}) made such a remarkable use of this mode that there is nothing to add to that. Therefore we shall carefully avoid using it -- unless it is hidden in a superimposition of modes that will make it unrecognizable." (The translations from the French are mine.)}

\noindent {\it Mode 2.---} In the language of Messiaen, this
mode is ``three times transposable", which means that, as a subset of
$\mathbb{Z}_{12}$, it is invariant under the translation $x\mapsto x+3$.
It is given by  the sequence of notes 
 {\it C, D$\flat$, E$\flat$,
E$\natural$, F$\sharp$, G, A, B$\flat$},  which corresponds to the {\it pc-set}
[0, 1, 3, 4, 6, 7, 9, 10], that is, the union of the orbits of $0$ and
$1$ by the maps $x\mapsto x+3$. It consists of 4 groups of 2 consecutive
elements in $\mathbb{Z}_{12}$. Messiaen writes, in \cite{Messiaen1944} p. 13 that ``modes of limited transposition are divisible in small symmetric groups". He points out that this mode has been
used by Stravinsky at the state of a rough sketch.  
In musical notation, this mode corresponds to the following musical representation:
\[
\includegraphics[width=.8\linewidth]{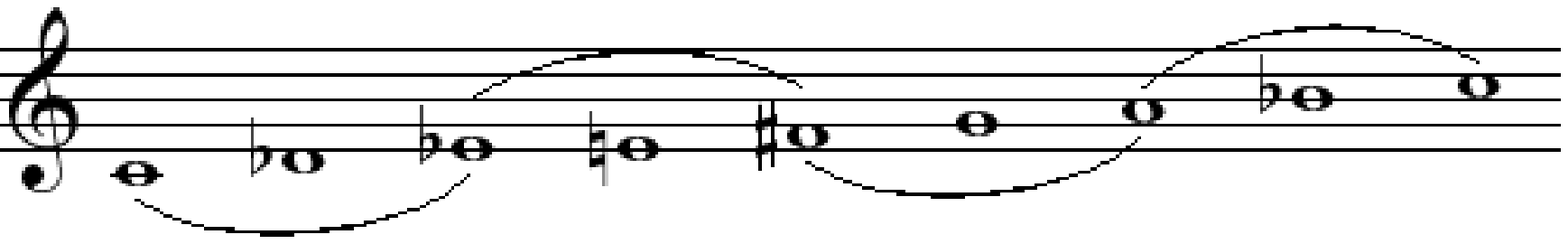}
\]
Messiaen used the second mode in his early compositions, for instance in the \emph{Banquet C\'eleste} (1928). This mode, and the subsequent ones, were extensively used by him.

\noindent {\it Mode 3.---} This mode is ``four times
transposable", and it is defined by the sequence of notes  {\it C, D, E$\flat$,
E$\natural$, F$\sharp$, G, A$\flat$, B$\flat$, B$\natural$}. 
In $\mathbb{Z}_{12}$, that is,
the sequence of integers mod 12 [0, 2, 3, 4, 6, 7, 8, 10, 11], 
which is the union of the
orbits of 0, 2 and 3 by the map $x\mapsto x+4$. The pc-set consists in 3 groups of
3 consecutive elements in $\mathbb{Z}_{12}$.
In musical notation this mode is described as follows:
\[
\includegraphics[width=.78\linewidth]{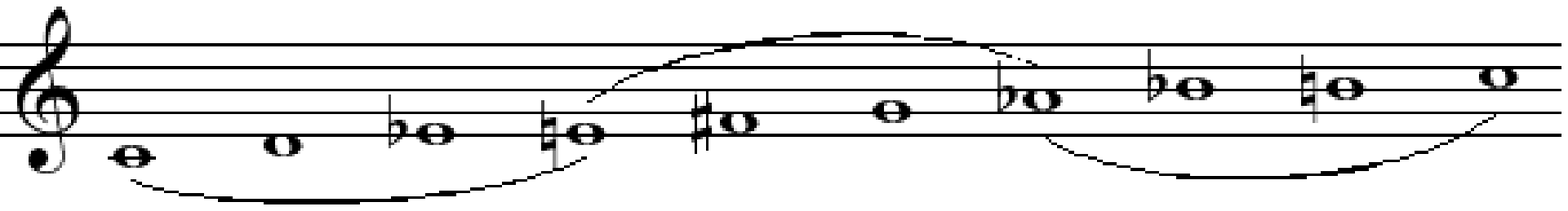}
\]

The remaining four modes are all ``six times transposable", which means that 
 as
subsets of the group $\mathbb{Z}_{12}$, they are invariant by the translation
$x \mapsto x+6$ mod 12. They are obtained  by taking in various ways unions
of  orbits of points by this map.

\medskip

\noindent {\it Mode 4.---} This is the sequence
{C, D$\flat$, D$\natural$, F, F$\sharp$, G, A$\flat$, B},
which corresponds to the {\it pc-set} 
[0, 1, 2, 5, 6, 7, 8, 11].
It is the union of the orbits of 0, 1, 2 and 5
by the map $x\mapsto x+6$.
This mode corresponds to the following intervals:
\[
\includegraphics[width=.8\linewidth]{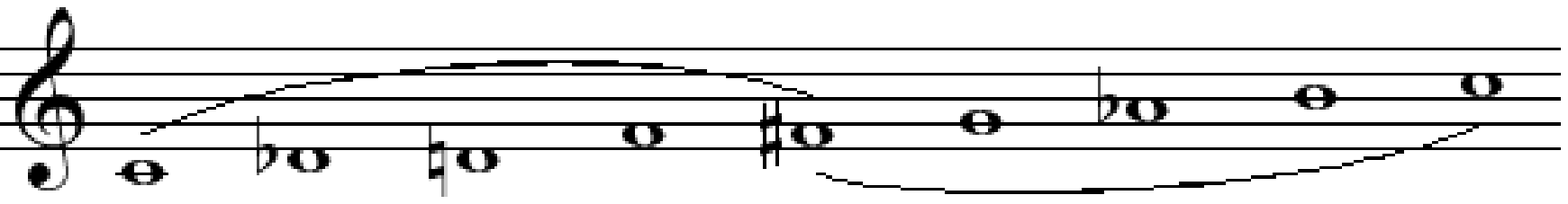}
\] 

\noindent {\it Mode 5.---}This is the sequence
{C, D$\flat$, F, F$\sharp$, G, B},
which corresponds to the {\it pc-set} 
[0, 1, 5, 6, 7, 11].
It is the union of the orbits of 
0, 1 and 5
by the map $x\mapsto x+6$.
In musical notation this mode is described as follows:
\[
\includegraphics[width=.72\linewidth]{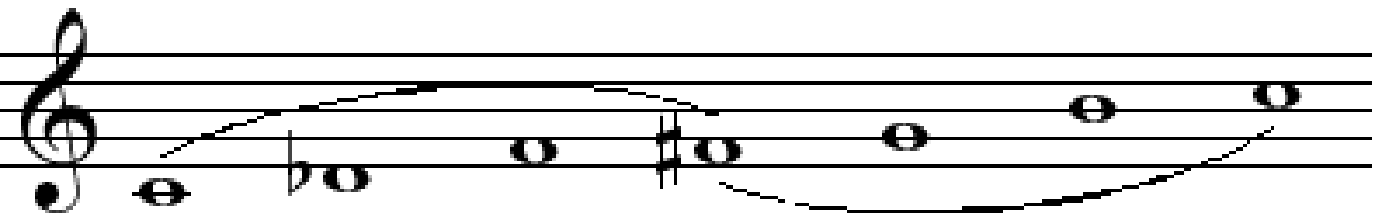}
\] 

\noindent {\it Mode 6.---} This is the sequence
{C, D, E, F, F$\sharp$, G$\sharp$, A$\sharp$, B}, 
which corresponds to the {\it pc-set}
[0, 2, 4, 5, 6,  8, 10, 11].
It is the union of the orbits of 0, 2, 4 and 5
by the map $x\mapsto x+6$.
In musical notation this mode is described by the following notes:
\[
\includegraphics[width=.8\linewidth]{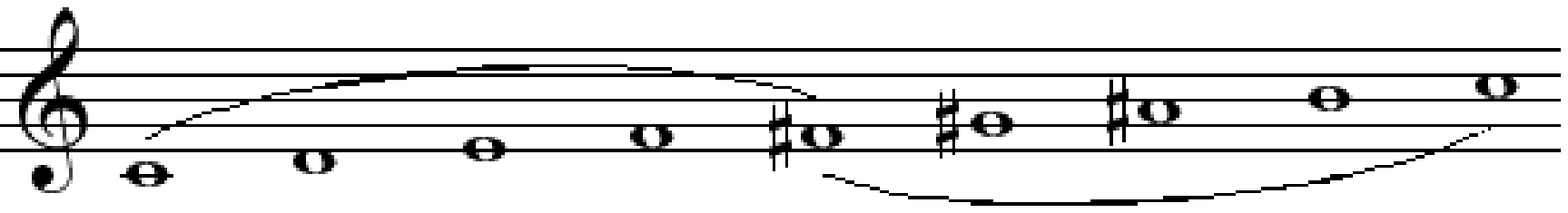}
\]

\noindent {\it Mode 7.---}
This is the sequence
{\it C, D$\flat$, D$\natural$, E$\flat$, E$\natural$, F$\sharp$, G, A$\flat$,
A$\natural$, B$\flat$},  which corresponds to the {\it pc-set}
[0, 1, 2, 3, 4, 6, 7, 8, 9, 10].
It is the union of the orbits of 0, 1, 2, 3 and 4
by the map $x\mapsto x+6$.
In musical notation this mode is described as follows:
\[
\includegraphics[width=.8\linewidth]{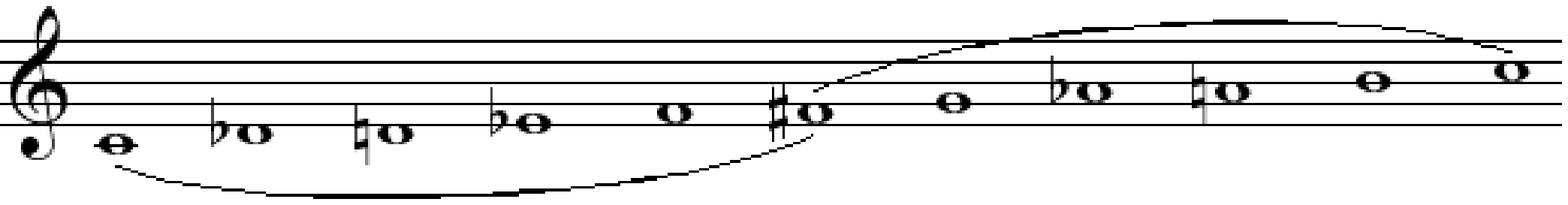}
\]

Each mode has its own transpositions, whose number is always less than twelve (by the definition of a mode of limited transposition, and unlike the usual major and minor scales).

The reader may notice that there are other divisions of the octave that produce modes of limited transposition, and Messiaen did not mention them. For instance, the complement in $\mathbb{Z}_{12}$ of some of the systems he defined is a mode of limited transposition which does not belong to this list of seven modes. Strangely enough, Messiaen declares in \cite{Messiaen1944}, vol. I, p. 51, after he presents his seven modes: ``The series is closed. It is mathematically impossible to find others, at least in our tempered system of 12 semi-tones." This is either due to the fact that the modes that are left over have a small number of notes, and they can be considered as truncations of the seven modes that he highlights (but it seems that Messiaen never mentioned this fact), or Messiaen simply missed them. We recall that Messiaen, despite the fact that he liked mathematics, had a very limited mathematical background. He also had no contact with mathematicians, even though some of his students, like Pierre Boulez and Iannis Xenakis\footnote{Boulez had one year of college mathematics and Xenakis was an engineer.} had some sort of it; but these students came several years after he developed the theory of modes of limited transposition. Messiaen continues in \cite{Messiaen1944}, vol. I by mentioning an analogous series of modes in the tempered system of quarter tones: ``In the system of tempered quarter tones, advocated by Haba and Wyshnegradsky, there exists a corresponding series." In fact, these modes of limited transposition in the microtonal universe were used by Georgui Rimsky-Korsakov, and more recently by Alain Louvier. They are classified in Jedrzejewski's book \cite{Jedr}.

A careful analysis of the intervallic structure of Messiaen's modes of limited transposition is made by Andreatta in \cite{A-Calcul} (p. 453 ff.)

Let us summarize. The musical atmosphere of a piece based on a mode of limited transposition can be characterized by a certain harmonious distribution of the notes which contributes to an ambiguous or mysterious ambience, which has no tonality, due (unlike the case of the classical compositions described by major and minor scales) to the lack of central notes, called classically the \emph{tonic} and the \emph{root}. 
 This is a consequence of the symmetry inherent in the scale used. The piece becomes aerial, like an impressionistic painting. In fact, the modes appear in Messiaen's works at the level at harmonies (that is, as notes played at the same time) rather than at the level of melodies. In his dialogues with Pierre Samuel  \cite{Messiaen1986}, Messiaen says:``The classical scale had a tonic. The antique modes had a finalis. My modes have neither a tonic nor a finalis. They are colors." In a famous talk in Kyoto, referred to as the ``Kyoto lecture" (``Conf\'erence de Kyoto") \cite{Messiaen-Kyoto}, Messiaen says: ``Modes are colored spots, like small colored countries, where the general color stays the same as long as we do not change the mode and the transposition." And also: ``People have often quoted my modes of limited transposition as scales. They are not scales, but harmonic colors".

The Japanese composer T\=oru Takemitsu  (1930-1996), the French composer Mich\`ele Reverdy (born in 1943), the British composer Michael Tippett (1905-1998) and several others  used Messiaen's modes in several of their pieces.

\subsection{\bf Modes and colors}

 We already quoted Messiaen talking about
the {\it voluptuous charm} produced by each of his modes, a charm,
which resides, as he says, in the {\it mathematical impossibility} of transposing the mode more than a small number of times (there are at most six transpositions). Beyond  the mathematical formalism, there is music,
and color. As we already alluded to, Messiaen had the rare faculty of having an inner vision of sounds as colors, and mixtures of sounds with mixtures of colors. According to this correspondence,
each of the  seven modes, at each of its transpositions,
has  a definite color.
Messiaen not only was capable of seeing colors while listening to a piece of music, but also when reading a music score,  that is, while hearing mentally the sounds (see \cite{Messiaen-Traite}, Tome 3, p. 79).
For instance, Mode  $2^3$ (that is, the second mode at its minor-second transposition, which
corresponds to the sequence of notes {\it C, D, E$\flat$, F, F$\sharp$,
A$\flat$, A$\natural$, B}), is a ``green
light foliage  with blue, silver and reddish
spots".  The colors that he perceived were, like music, 
changing  and moving. Like musical lines, these colors are superimposed   and interlaced.
In the Preface to his {\it Couleurs de la Cit\'e C\'eleste}, he
writes about this piece: 
\begin{quote}\small 
The melodic or rhythmic themes,
the complexes of sounds and of timbres, evolve in the same way as colors do. In
their perpetually renewed variations, we can find (by analogy) warm and cold
tints, complementary colors influencing their neighbors,  
with a gradation towards  white, turned down by black, etc.
\end{quote}

The intuitions on the relation sound-color are present at
different levels in the works various artists. One can recall here  
  Baudelaire's  {\it Correspondences},   described in the
  following quatrain derived from a poem of {\it Les Fleurs du Mal}:
\begin{quote}\small
{\it Comme de longs \'echos qui de loin se confondent
\\
Dans une t\'en\'ebreuse et profonde unit\'e
\\
Vaste comme la nuit et comme la clart\'e,
\\
Les parfums, les couleurs et les sons se r\'epondent.}
\end{quote}

\medskip
In 
his article {\it R. Wagner et Tannh\"auser \`a Paris} (1861), Baudelaire says that 
``it would have been very surprising if sound could not suggest color, if
color could not give the idea of a melody, etc."

One can  also mention the mystical 
{\it Illuminations} of
Rimbaud, which he describes in  {\it Une Saison en
Enfer}: ``J'inventai la couleur des voyelles! A noir, E blanc, I rouge, O bleu, U
vert,..."

 Sch\"onberg, in his piece {\it Die gl\"uckliche Hand}
 (1910-1913), gave specific instructions on colors
 to be projected on  a screen, during its public performance. 
Likewise, Scriabin required that light beams 
of various colors accompany
 his symphonic poem   {\it Prometh\'ee, le po\`eme de feu} (1908-1910), operated by a  device connected to various keys of the organ. 
 Back in the Century of Enlightenment, the mathematician Castel whom we already mentioned
 invented a ``color organ" ({\it clavecin oculaire}), which he describes 
 in detail in his
 paper \cite{Castel}.
We also remember now that Newton analyzed the relation of music with color 
by comparing the spectrum of sound to the
spectrum of light, and we can even trace back this 
relation to Aristotle's {\it De sensu}. Thus, we are back to the origins. In fact,
 the physical perception of
musical objects (tones, modes, chords, etc.) as colors remains an interesting subject of
research, from the
physiological  and from the
psychological points of view,
and mathematics is, as usual, at the basis.


\end{document}